\documentclass[10pt]{article}

\usepackage{authblk}
\usepackage[text={6.75in,9.5in},centering,letterpaper]{geometry}
\usepackage{amsfonts,amsmath}
\usepackage{paralist}
\usepackage{graphicx,caption}
\usepackage[numbers,comma,square,sort&compress]{natbib}

\setcaptionmargin{0.25in}

\newtheorem{Lemma}{Lemma}[section]
\newtheorem{Hypothesis}{Hypothesis}
\newtheorem{Remark}{Remark}
\newtheorem{Theorem}{Theorem}

\newenvironment{Proof}[1][.]%
 {\begin{trivlist}\item[]\textbf{Proof#1 }}%
 {\hspace*{\fill}$\rule{0.3\baselineskip}{0.35\baselineskip}$\end{trivlist}}

\newenvironment{Acknowledgment}%
 {\begin{trivlist}\item[]\textbf{Acknowledgments }}{\end{trivlist}}

\numberwithin{equation}{section}

\def\Re{\mathop\mathrm{Re}\nolimits}
\def\Im{\mathop\mathrm{Im}\nolimits}
\def\sech{\mathop\mathrm{sech}\nolimits}
\def\Rg{\mathop\mathrm{Rg}\nolimits}

\def\id{\mathbb{I}}

\newcommand{\C}{\mathbb{C}}
\newcommand{\R}{\mathbb{R}}
\newcommand{\rme}{\mathrm{e}}
\newcommand{\rmi}{\mathrm{i}}
\newcommand{\rmO}{\mathrm{O}}

\title{Computing Evans functions numerically via boundary-value problems}

\author[1]{Blake Barker}
\author[2]{Rose Nguyen}
\author[3]{Bj\"orn Sandstede}
\author[4]{Nathaniel Ventura}
\author[5]{Colin Wahl}

\affil[1]{Department of Mathematics, Brigham Young University, Provo, UT 84604, USA}
\affil[2]{Department of Applied Mathematics, University of Washington, Seattle, WA 98195, USA}
\affil[3]{Division of Applied Mathematics, Brown University, Providence, RI 02912, USA}
\affil[4]{Department of Mathematics, Binghamton University, Binghamton, NY 13902, USA}
\affil[5]{Applied Science and Technology, University of California at Berkeley, Berkeley, CA 94720, USA}

\begin{document}
\maketitle
	
\begin{abstract}
The Evans function has been used extensively to study spectral stability of travelling-wave solutions in spatially extended partial differential equations. To compute Evans functions numerically, several shooting methods have been developed. In this paper, an alternative scheme for the numerical computation of Evans functions is presented that relies on an appropriate boundary-value problem formulation. Convergence of the algorithm is proved, and several examples, including the computation of eigenvalues for a multi-dimensional problem, are given. The main advantage of the scheme proposed here compared with earlier methods is that the scheme is linear and scalable to large problems.
\end{abstract}

% --------------------------------------------------------------------------------------

\section{Introduction}\label{s1}

In this paper, we explore ways to study the spectral stability of travelling waves in spatially extended partial differential equations (PDEs). To provide background, we consider reaction-diffusion systems
\begin{equation}\label{e1}
u_t = D u_{xx} + f(u), \quad x\in\R, \quad u\in\R^n,
\end{equation}
though we emphasize that the methods presented below are applicable to a much larger class of PDEs. Travelling waves are solutions of (\ref{e1}) of the form $u(x,t)=u_*(x-ct)$, where $c$ is the wave speed associated with the profile $u_*(x)$. For simplicity, we focus in the introduction on pulses for which $\lim_{|x|\to\infty}u_*(x)=0$. Linearizing (\ref{e1}) about the pulse profile $u_*(x)$ in a frame that moves with the travelling pulse gives the operator
\[
\mathcal{L}u := D u_{xx} + c u_x + f_u(u_*(x)) u
\]
posed, for instance, on $L^2(\R,\R^n)$. We say that the travelling pulse is spectrally stable if the spectrum of $\mathcal{L}$ is contained entirely in the open left half-plane with the exception of a simple eigenvalue at the origin that is enforced by translational symmetry. It has been shown, for instance in \citep{He}, that spectral stability implies nonlinear stability of the pulse with respect to (\ref{e1}) under sufficiently small perturbations. It is often difficult to prove spectral stability, and numerical computations of the spectrum of $L$ are therefore frequently the only way to determine stability of a given travelling wave. A natural approach to computing the spectrum of $\mathcal{L}$ numerically is to replace the unbounded domain $\R$ by a large but finite interval, add boundary conditions to make the resulting problem well-posed, discretize the operator using finite difference or spectral schemes, and apply eigenvalue problem solvers to the resulting large matrix. In many cases, proceeding in this fashion will produce reliable and accurate results. There are, however, several caveats to consider. Firstly, the limit of the spectra obtained on intervals of the form $(-\ell,\ell)$ as $\ell\to\infty$ may not coincide with the spectrum of $\mathcal{L}$ posed on $\R$; see \citep{SS2000}. A related issue is the presence of pseudo spectra that can lead to spurious eigenvalues when using iterative solvers. Secondly, one is often interested in identifying situations where eigenvalues can emerge from the essential spectrum \citep{KS1998,KS2002}: in this case, eigenvalue problem solvers do not help. For these reasons, we follow a different approach in this paper that focuses on the Evans function.

To review the approach via Evans functions, we write the eigenvalue problem $\mathcal{L}u=\lambda u$ associated with the operator $\mathcal{L}$ as the linear ordinary differential equation
\[
\begin{pmatrix} u\\ u_x \end{pmatrix}_x = \begin{pmatrix} 0 & 1 \\ D^{-1}(\lambda-f_u(u_*(x))) & -cD^{-1} \end{pmatrix}
\begin{pmatrix} u\\ u_x \end{pmatrix},
\]
which we will write as
\begin{equation}\label{e2}
U_x = A(x,\lambda) U, \qquad U\in\C^{2n}.
\end{equation}
Note that $A(x,\lambda)$ will converge to a matrix $A_\infty(\lambda)$ as $x\to\pm\infty$, since we assumed that the profile converges to zero as $|x|\to\infty$. As long as $\lambda\in\C$ is not in the essential spectrum of $\mathcal{L}$, the matrix $A_\infty(\lambda)$ will be hyperbolic, that is, does not have eigenvalues on the imaginary axis, and nontrivial bounded solutions to (\ref{e2}) will therefore automatically decay exponentially as $|x|\to\infty$. To find eigenvalues of $\mathcal{L}$ away from the essential spectrum, we can therefore define the spaces
\[
E^\pm(\lambda) := \{ U_0\in\C^{2n};\; U(x) \mbox{ satisfies (\ref{e2}) with } U(0)=U_0 \mbox{, and } U(x)\to0  \mbox{ as }x\to\pm\infty \}
\]
and define the Evans function $E(\lambda)$ via $E(\lambda)=E^-(\lambda)\wedge E^+(\lambda)$, where $\wedge$ denotes the wedge product of two vector spaces: alternatively, we can interpret the Evans function as the Wronskian associated with the set of solutions of \eqref{e2} that decay as $x\to\infty$ and the set of solutions that decay as $x\to-\infty$. The Evans function $E(\lambda)$ is analytic in $\lambda$ and vanishes precisely at eigenvalues of $\mathcal{L}$. As shown in \citep{AGJ}, the multiplicity of a root $\lambda$ corresponds to the multiplicity of $\lambda$ as a PDE eigenvalue of $\mathcal{L}$. In particular, we can determine the number of eigenvalues in a given region in the complex plane by computing an appropriate winding number of $E(\lambda)$.

Computing the Evans function numerically amounts to finding approximate basis vectors of $E^\pm(\lambda)$ that vary analytically in $\lambda$. A key challenge is that solutions integrated forward or backward will collapse on the dominant growing modes, which makes it difficult to compute bases. The exterior-product or compound-matrix method, used as early as 1995 \citep{AS} and later developed further in \citep{Allen2002,BDG,Br1,Br2,BrZ}, addresses this challenge by lifting (\ref{e2}) to the Grassmannian, that is the manifold of vector spaces, on which $E^\pm(\lambda)$ correspond to a single maximally unstable or stable mode. Many systems have been successfully studied using this approach, but the key limitation is that the dimension of the resulting system on the Grassmannian is typically of size $\binom{2n}{n}$, which is too large in practice to be used for systems with $n$ larger than two or three. In 2006, this dimensionality challenge was addressed by finding basis vectors of $E^\pm(\lambda)$ via solving (\ref{e2}) using continuous orthogonalization whereby an orthonormal basis for the desired manifold is evolved along with the determinant of the coordinate matrix, thus allowing recovery of an analytic Evans function via a numerically stable algorithm \citep{HuZ2,Z5}. Since then, this method has been used to study a range of systems; see \citep{BFZ, BHLZ,BJNRZ0,BaLeZ,HLyZ1} for a few examples. Other low-dimensional shooting approaches followed, such as schemes that utilize the relation between the Grassmann and Stiefel manifolds \citep{Ledoux2009,Ledoux2010}.

Whether variants of the compound-matrix method or continuous orthogonalization are used, the resulting systems are effectively nonlinear, whilst the computation of the ingredients of the Evans function is linear. Our goal in this paper is to propose an algorithm that is linear and scalable, and that can therefore be used for multi-dimensional problems, where $x$ lies in a cylindrical domain. We will achieve this by formulating the problem of finding basis vectors of $E^\pm(\lambda)$ as an appropriate linear boundary-value problem.

The remainder of this paper is organized as follows. In \S\ref{s2}, we will give two constructions of the Evans function and outline how these can be used to design stable and accurate numerical algorithms for the computation of the Evans function. In \S\ref{s3}, we will apply our algorithm to three test problems to demonstrate its accuracy and scalability. We conclude in \S\ref{s4} with a discussion of open problems.

% --------------------------------------------------------------------------------------

\section{Evans functions and their numerical computation}\label{s2}

In this section, we will describe the setting in which we will work throughout the remainder of this paper, recall the definition of the Evans function, provide an alternative formulation that will form the basis of the proposed numerical framework, and, finally, explain our algorithm and prove its convergence. We remark that our notation will differ from the one used in the introduction.

Throughout this paper, we consider the linear system
\begin{equation}\label{e10}
U_x = A(x,\lambda) U, \qquad U\in\C^n,
\end{equation}
and assume that the following hypothesis is met.

\begin{Hypothesis}\label{h1}
We assume that there is an open, bounded, and simply connected set $\Omega\subset\C$, a continuous matrix-valued function $A: \R\times\bar\Omega\to\C^{n\times n}, (x,\lambda)\mapsto A(x,\lambda)$, and continuous functions $A_\pm:\bar\Omega\to\C^{n\times n}, \lambda\mapsto A_\pm(\lambda)$ with the following properties:
\begin{compactenum}
\item $A(x,\lambda)$ and $A_\pm(\lambda)$ are analytic in $\lambda$ for $\lambda\in\Omega$.
\item $A(x,\lambda)\to A_\pm(\lambda)$ exponentially as $x\to\pm\infty$.
\item $A_\pm(\lambda)$ are hyperbolic for $\lambda\in\bar\Omega$ and have precisely $r$ eigenvalues, counted with multiplicity, with negative real part.
\end{compactenum}
\end{Hypothesis}

We will also require the following assumption for some of our results as it simplifies our arguments, though we stress that our algorithm can be modified to hold even if this hypothesis is not met.

\begin{Hypothesis}\label{h2}
The eigenvalues of the matrices $A_\pm(\lambda)$ appearing in \ref{h1} are simple for all $\lambda\in\bar\Omega$.
\end{Hypothesis}

We denote the spectral projections associated with the eigenvalues of $A_\pm(\lambda)$ that have negative real part by $P^s_\pm(\lambda)$ and write $P^u_\pm(\lambda)$ for the complementary spectral projections onto the generalized eigenspaces associated with unstable eigenvalues. Note that these projections are analytic in $\lambda$ for $\lambda\in\Omega$ and that $\dim \Rg(P^s_\pm(\lambda))=r$ and therefore $\dim \Rg(P^u_\pm(\lambda))=n-r$.

% --------------------------------------------------------------------------------------

\subsection{The Evans function}\label{s2.1}

We assume \ref{h1}. It follows from this hypothesis that the space
\[
E^s_+(y,\lambda) := \{ U_0\in\C^n:\; U(x) \mbox{ satisfies (\ref{e10}) with } U(y)=U_0 \mbox{, and } U(x)\to0  \mbox{ as }x\to\infty \}
\]
has dimension $r$, is analytic in $\lambda\in\Omega$, and satisfies $E^s_+(x,\lambda)\to\Rg(P^s_+(\lambda))$ as $x\to\infty$ for each $\lambda\in\bar\Omega$. In particular, there is a number $L_0\gg1$ that does not depend on $\lambda\in\bar\Omega$ such that the projection onto $E^s_+(L,\lambda)$ along $\Rg(P^u_+(\lambda))$ is well defined and analytic in $\lambda\in\Omega$ for each $L\geq L_0$. Using this projection and applying the results in \citep[Ch.~II.4.2]{Kato} to it, we can construct a basis $\{Z_j^+(L,\lambda)\}_{j=1,\ldots,r}$ of $E^s_+(L,\lambda)$ that depends analytically on $\lambda\in\Omega$. We denote by $Z_j^+(x,\lambda)$ the solutions of (\ref{e10}) with initial conditions $Z_j^+(L,\lambda)$ at $x=L$.

The same step can be applied to (\ref{e10}) for $x\leq0$ resulting in an analytically varying basis $\{Z_j^-(-L,\lambda)\}_{j=1,\ldots,n-r}$ of $E^u_-(-L,\lambda)$ where
\[
E^u_-(y,\lambda) := \{ U_0\in\C^n:\; U(x) \mbox{ satisfies (\ref{e10}) with } U(y)=U_0 \mbox{, and } U(x)\to0  \mbox{ as }x\to-\infty \}.
\]
The Evans function is now defined as
\begin{equation}\label{e21}
E_L(\lambda) := \det(Z_1^-(0,\lambda),\ldots,Z_{n-r}^-(0,\lambda),Z_1^+(0,\lambda),\ldots,Z_r^+(0,\lambda))
\rme^{-\mathrm{tr}(A^-(\lambda)P^u_-(\lambda))L} \rme^{\mathrm{tr}(A^+(\lambda)P^s_+(\lambda))L},
\end{equation}
where the exponential terms on the right-hand side ensure that the Evans function stays bounded as $L$ increases. Since (\ref{e10}) and the initial data $Z_j^\pm(\pm L,\lambda)$ depend analytically on $\lambda$ for $\lambda\in\Omega$, the Evans function $E_L:\Omega\to\C$ is analytic. We emphasize that the roots of $E_L(\lambda)$ do not depend on $L$.

% --------------------------------------------------------------------------------------

\subsection{A boundary-value problem formulation of the Evans function}\label{s2.2}

We now provide a different construction of the Evans function: we assume that \ref{h1} and \ref{h2} are met. We focus initially solely on existence for $\lambda\in\bar\Omega$, rather than analyticity, and will show in \S\ref{s2.3} how analyticity is recovered.

First, consider solutions of (\ref{e10}) on $\R^+$. We order the $n$ distinct eigenvalues $\nu_j^+(\lambda)$ of $A^+(\lambda)$ by their real part so that $\Re\nu^+_j(\lambda)<0$ for $j=1,\ldots,r$. We can also assume that these eigenvalues vary continuously in $\lambda\in\bar\Omega$ and then choose associated eigenvectors $V^+_j(\lambda)$ of $A^+(\lambda)$ so that these also vary continuously in $\lambda\in\bar\Omega$ for $j=1,\ldots,r$. We denote the identity matrix by $\id$.

\begin{Lemma}[{\citep[Ch.~3.8]{CodLev}}]\label{l1}
Assume that \ref{h1} and \ref{h2} are met. For each $j=1,\ldots,r$, the system
\begin{equation}\label{e11}
V_x = (A(x,\lambda)-\nu^+_j(\lambda)\id) V
\end{equation}
has a solution $V^+_j(x,\lambda)$ on $\R^+$ that varies continuously in $\lambda\in\bar\Omega$ and satisfies $V^+_j(x,\lambda)\to V^+_j(\lambda)$ as $x\to\infty$ for each $\lambda\in\bar\Omega$.
\end{Lemma}

Note that, since $V^+_j(x,\lambda)$ satisfies (\ref{e11}) for $j=1,\ldots,r$, we know that each function
\begin{equation}\label{e13}
U(x) = \rme^{\nu^+_j(\lambda)x} V^+_j(x,\lambda), \qquad j=1,\ldots,r
\end{equation}
satisfies (\ref{e10}) with $U(x)\to0$ as $x\to\infty$. Furthermore, the vectors $V^+_j(0,\lambda)$ with $j=1,\ldots,r$ are linearly independent and continuous in $\lambda\in\bar\Omega$. Combining these results shows that the vectors $V^+_j(x,\lambda)$ with $j=1,\ldots,r$ form a continuously varying basis of $E^s_+(x,\lambda)$ for $\lambda\in\bar\Omega$ and each $x\geq0$.

Proceeding in the same way for $x\in\R^-$, where we order the eigenvalues $\nu^-_j(\lambda)$ according to $\Re\nu^-_j(\lambda)>0$ for $j=1,\ldots,n-r$ and $\Re\nu^-_j(\lambda)<0$ otherwise, we arrive at a continuously varying basis $\{V^-_j(0,\lambda)\}_{j=1,\ldots,n-r}$ of $E^u_-(0,\lambda)$. We can now define our alternative Evans function $\mathcal{E}(\lambda)$ via
\begin{equation}\label{e22}
\mathcal{E}(\lambda) := \det(V_1^-(0,\lambda),\ldots,V_{n-r}^-(0,\lambda),V_1^+(0,\lambda),\ldots,V_r^+(0,\lambda))
\end{equation}
that is defined and continuous for $\lambda\in\bar\Omega$. Note that $\mathcal{E}(\lambda)$ does not depend on any additional parameters and that it vanishes if, and only if, $E_L(\lambda)$ vanishes. However, we do not know whether $\mathcal{E}(\lambda)$ is analytic.

% --------------------------------------------------------------------------------------

\subsection{Comparison of the two formulations}\label{s2.3}

We focus again first on solutions defined on $\R^+$. Let $Z^+(x,\lambda):=(Z^+_1(x,\lambda),\ldots,Z^+_r(x,\lambda))\in\C^{n\times r}$ be the matrix with columns given by the solutions $Z^+_j(x,\lambda)$ discussed in \S\ref{s2.1} for $j=1,\ldots,r$ and similarly define $V^+(x,\lambda)$ to be the matrix with columns consisting of the solutions $V^+_j(x,\lambda)$ discussed in \S\ref{s2.2} for $j=1,\ldots,r$. Since both sets of column vectors form a basis of $E^s_+(x,\lambda)$ for all $\lambda\in\Omega$ and $x\geq0$, we know that there is a unique matrix $C^+_L(\lambda)\in\C^{r\times r}$ so that
\begin{equation}\label{e12}
V^+(L,\lambda) C^+_L(\lambda) = Z^+(L,\lambda)
\end{equation}
for all $\lambda\in\Omega$. Applying the fundamental matrix solution $\Phi(x,L,\lambda)$ of (\ref{e10}) to both sides of (\ref{e12}), we see that
\begin{equation}\label{n1}
\Phi(0,L,\lambda) V^+(L,\lambda) C^+_L(\lambda) = Z^+(0,\lambda).
\end{equation}
Using the relationship (\ref{e13}) between solutions of (\ref{e10}) and (\ref{e11}), we see that
\[
\Phi(x,L,\lambda) \rme^{\nu^+_j(\lambda)L} V^+_j(L,\lambda) = \rme^{\nu^+_j(\lambda)x} V^+_j(x,\lambda)
\]
and therefore
\[
\Phi(0,L,\lambda)V^+_j(L,\lambda) = \rme^{-\nu^+_j(\lambda)L} V^+_j(0,\lambda), \qquad j=1,\ldots,r.
\]
Thus, denoting by $D^+_L(\lambda)\in\C^{r\times r}$ the diagonal matrix with entries $\rme^{-\nu^+_j(\lambda)L}$ on the diagonal, we have
\begin{equation}\label{e62}
\Phi(0,L,\lambda) V^+(L,\lambda) = V^+(0,\lambda) D^+_L(\lambda)
\end{equation}
for all $\lambda$. Note that
\begin{equation}\label{n2}
\det D^+_L(\lambda) = \rme^{-\mathrm{tr}(A^+(\lambda)P^s_+(\lambda))L}.
\end{equation}
Proceeding analogously for the matrices formed by the solutions defined for $x\leq0$, we therefore have
\begin{eqnarray*}
E_L(\lambda) & = & \det(Z^-(0,\lambda),Z^+(0,\lambda))
\rme^{-\mathrm{tr}(A^-(\lambda)P^u_-(\lambda))L} \rme^{\mathrm{tr}(A^+(\lambda)P^s_+(\lambda))L} \\ & \stackrel{(\ref{n1})}{=} &
\det(\Phi(0,-L,\lambda)V^-(-L,\lambda)C^-_L(\lambda),\Phi(0,L,\lambda)V^+(L,\lambda)C^+_L(\lambda))
\rme^{-\mathrm{tr}(A^-(\lambda)P^u_-(\lambda))L} \rme^{\mathrm{tr}(A^+(\lambda)P^s_+(\lambda))L} \\ & \stackrel{(\ref{e62})}{=}  &
\det(V^-(0,\lambda)D^-_L(\lambda)C^-_L(\lambda),V^+(0,\lambda)D^+_L(\lambda)C^+_L(\lambda))
\rme^{-\mathrm{tr}(A^-(\lambda)P^u_-(\lambda))L} \rme^{\mathrm{tr}(A^+(\lambda)P^s_+(\lambda))L} \\ & = &
\det\left[ (V^-(0,\lambda),V^+(0,\lambda)) 
\begin{pmatrix} D^-_L(\lambda)C^-_L(\lambda) & 0 \\ 0 & D^+_L(\lambda)C^+_L(\lambda) \end{pmatrix}\right]
\rme^{-\mathrm{tr}(A^-(\lambda)P^u_-(\lambda))L} \rme^{\mathrm{tr}(A^+(\lambda)P^s_+(\lambda))L} \\ & = &
\det(V^-(0,\lambda),V^+(0,\lambda))
\det(D^+_L(\lambda)) \det(D^-_L(\lambda))
\det(C^+_L(\lambda)) \det(C^-_L(\lambda)) \\ && \times
\rme^{-\mathrm{tr}(A^-(\lambda)P^u_-(\lambda))L} \rme^{\mathrm{tr}(A^+(\lambda)P^s_+(\lambda))L} \\ & \stackrel{(\ref{n2})}{=} &
\underbrace{\det(V^-(0,\lambda),V^+(0,\lambda))}_{=:\mathcal{E}(\lambda)} \underbrace{\det(C^+_L(\lambda)) \det(C^-_L(\lambda))}_{=:\mathcal{C}_L(\lambda)} \\ & = &
\mathcal{E}(\lambda) \mathcal{C}_L(\lambda).
\end{eqnarray*}
Thus, the analytic function $E_L(\lambda)$ can be calculated from $\mathcal{E}(\lambda)$ provided we can determine $C^\pm_L(\lambda)$ from (\ref{e12}). Note that $\mathcal{E}(\lambda)$ does not depend on $L$ and that there are constants $C_\pm>0$ such that
\begin{equation}\label{e71}
0 < C_- \leq |\det C^\pm_L(\lambda)| \leq C_+
\end{equation}
for all $L\geq L_0$ and $\lambda\in\bar\Omega$ by construction. We will now use these results to design our numerical algorithm.

% --------------------------------------------------------------------------------------

\subsection{Numerical algorithm}\label{s2.4}

Assuming again \ref{h1} and \ref{h2}, we can now outline our numerical algorithm for the computation of the Evans function $E_L(\lambda)$ via the numerical approximation of $\mathcal{E}(\lambda)$ and $\mathcal{C}_L(\lambda)$ introduced in the last two sections. We will first describe the algorithm and then present a theorem that states that the numerical Evans function is analytic and approximates the exact Evans function.

Throughout the remainder, we denote by ${ }^*$ the complex-conjugate transpose of a number, vector, or matrix, and denote by $\langle U,V\rangle:=U^* \cdot V$ the scalar product in $\C^n$. Furthermore, we choose eigenvectors $V^+_j(\lambda)$ and $W^+_j(\lambda)$ belonging to the simple eigenvalues $\nu^+_j(\lambda)$ of $A^+(\lambda)$ and $\nu^+_j(\lambda)^*$ of $A^+(\lambda)^*$, respectively, that vary continuously in $\lambda\in\bar\Omega$. Finally, we fix a number $L\gg1$.

\paragraph{Step~1.} Use the spectral projection $P^s_+(\lambda)$ and Kato's algorithm \cite[Ch.~II.4.2]{Kato} to compute an analytic basis $\{Z^+_j(\lambda)\}_{j=1,\ldots,r}$  of the stable eigenspace $E^s_+(\lambda)$ of $A^+(\lambda)$.

\paragraph{Step~2.} Fix $\lambda\in\bar\Omega$, order the eigenvalues of $A^+(\lambda)$ according to their real part, and, starting with $j=1$ and ending at $j=r$, iteratively calculate solutions $V(x)=V^+_j(x,\lambda)$ of the linear system
\begin{compactenum}
\item $V_x=(A(x,\lambda)-\nu^+_j(\lambda)\id)V$ on $0<x<L$;
\item $\langle V^+_k(0,\lambda),V(0)\rangle=0$ for $k=1,\ldots,j-1$;
\item $\langle W^+_j(\lambda),V(L)\rangle=\langle W^+_j(\lambda),V^+_j(\lambda)\rangle$;
\item $\langle W^+_k(\lambda),V(L)\rangle=0$ for $k=j+1,\ldots,n$.
\end{compactenum}

\paragraph{Step~3.} Find $\tilde{C}^+_L(\lambda)\in\C^{r\times r}$ so that
\begin{equation}\label{e81}
V^+_L(L,\lambda) \tilde{C}^+_L(\lambda) = Z^+(\lambda)
\end{equation}
for all $\lambda\in\bar\Omega$, where $V^+_L(L,\lambda)$ and $Z^+(\lambda)$ are the $n\times r$ matrices with columns given by $V^+_j(L,\lambda)$ and $Z^+_j(\lambda)$, respectively.

\paragraph{Step~4.} Repeat steps~1-3 for $-L<x<0$ and $A^-(\lambda)$.

We then set
\begin{equation}\label{23}
\tilde{E}_L(\lambda) := \underbrace{\det(V^-_L(0,\lambda),V^+_L(0,\lambda))}_{=:\tilde{\mathcal{E}}_L(\lambda)} \underbrace{\det(\tilde{C}^-_L(\lambda)) \det(\tilde{C}^+_L(\lambda))}_{=:\tilde{\mathcal{C}}_L(\lambda)}.
\end{equation}

We have the following theorem that we will prove in \S\ref{s2.5}.

\begin{Theorem}\label{t1}
Assume \ref{h1} and \ref{h2}. There are constants $\eta,C,L_0>0$ so that the system described in steps~1-4 has a unique solution for each $L\geq L_0$ and
\begin{equation}\label{et1}
|\mathcal{E}(\lambda) - \tilde{\mathcal{E}}_L(\lambda)| + |\mathcal{C}_L(\lambda) - \tilde{\mathcal{C}}_L(\lambda)|\leq C \rme^{-\eta L}
\end{equation}
and therefore also
\begin{equation}\label{et2}
|E_L(\lambda) - \tilde{E}_L(\lambda)| \leq C \rme^{-\eta L}
\end{equation}
for $L\geq L_0$ uniformly in $\lambda\in\bar\Omega$. Furthermore, $\tilde{E}_L(\lambda)$ is analytic in $\lambda$ for $\lambda\in\Omega$ for each $L\geq L_0$.
\end{Theorem}

Note that the roots of $E_L(\lambda)$ are given by the roots of $\mathcal{E}(\lambda)$, which does not depend on $L$. Using (\ref{e71}) together with (\ref{et1}), we therefore conclude from Rouch\'{e}'s theorem that the roots of $E_L(\lambda)$ and $\tilde{E}_L(\lambda)$, counted with multiplicity, are $\rmO(\rme^{-\eta L})$ close to each other.

We can now use the numerical algorithm described in steps~1-4 to compute winding numbers of Evans functions. Choose a closed smooth curve $\Gamma$ in $\Omega$ that has no self-intersections and parametrize the curve by $\lambda(s)$ with $s\in[0,1]$. We can then carry out the calculations in steps~1-4 through numerical continuation in the parameter $s$, for instance by selecting a finite discretization $\{s_m\}_{m=1,\ldots,M}$ of $[0,1]$. Since our boundary-value problem is linear, we can also use parallel implementations to speed up the computations.

% --------------------------------------------------------------------------------------

\subsection{Proof of Theorem~\ref{t1}}\label{s2.5}

In this section, we prove Theorem~\ref{t1}: we assume that Hypotheses \ref{h1} and \ref{h2} are met.

First note that parts~(iii)-(iv) of step~2 of the algorithm described in \S\ref{s2.4} show that the vectors $V^+_j(L,\lambda)$ are linearly independent for $j=1,\ldots,r$ and lie in $E^s_+(\lambda)$. In particular, the matrix $\tilde{C}^+_L(\lambda)$ is well defined for all $\lambda\in\bar\Omega$. Next, we show that the numerical Evans function $E_L(\lambda)$ is analytic in $\lambda\in\Omega$. 

\begin{Proof}[ of analyticity of $\tilde{E}_L(\lambda)$.]
Recalling that $\Phi(x,y,\lambda)$ is the fundamental matrix solution of (\ref{e10}) and that we defined $D^+_L(\lambda)\in\C^{r\times r}$ to be the diagonal matrix with entries $\rme^{-\nu^+_j(\lambda)L}$ on the diagonal, we conclude from (\ref{e62}) that
\begin{equation}\label{n3}
\Phi(0,L,\lambda) V^+_L(L,\lambda) = V^+_L(0,\lambda) D^+_L(\lambda)
\end{equation}
for all $\lambda$, and analogous expressions for $V^-_L(-L,\lambda)$. Hence, 
\begin{eqnarray*}
\lefteqn{\det(\Phi(0,-L,\lambda)Z^-(\lambda),\Phi(0,L,\lambda)Z^+(\lambda)) } \\ & \stackrel{(\ref{e81})}{=} &
\det(\Phi(0,-L,\lambda)V^-_L(-L,\lambda)\tilde{C}^-_L(\lambda),\Phi(0,L,\lambda)V^+_L(L,\lambda)\tilde{C}^+_L(\lambda)) \\ & \stackrel{(\ref{n3})}{=} &
\det(V^-_L(0,\lambda)D^-_L(\lambda)\tilde{C}^-_L(\lambda),V^+_L(0,\lambda)D^+_L(\lambda)\tilde{C}^+_L(\lambda)) \\ &= &
\det\left[ (V^-_L(0,\lambda),V^+_L(0,\lambda)) \begin{pmatrix} D^-_L(\lambda) & 0 \\ 0 & D^+_L(\lambda) \end{pmatrix} \begin{pmatrix} \tilde{C}^-_L(\lambda) & 0 \\ 0 & \tilde{C}^+_L(\lambda) \end{pmatrix}\right] \\ & = &
\det(V^-_L(0,\lambda),V^+_L(0,\lambda)) \det(\tilde{C}^-_L(\lambda)) \det(\tilde{C}^+_L(\lambda))
\det(D^-_L(\lambda)) \det(D^+_L(\lambda)) \\ & \stackrel{(\ref{n2})}{=} &
\det(V^-_L(0,\lambda),V^+_L(0,\lambda)) \det(\tilde{C}^-_L(\lambda)) \det(\tilde{C}^+_L(\lambda))
\rme^{\mathrm{tr}(A^-(\lambda)P^u_-(\lambda))L} \rme^{-\mathrm{tr}(A^+(\lambda)P^s_+(\lambda))L}
\end{eqnarray*}
and therefore
\begin{eqnarray*}
\tilde{E}_L(\lambda) & = & \det(V^-_L(0,\lambda),V^+_L(0,\lambda)) \det(\tilde{C}^-_L(\lambda)) \det(\tilde{C}^+_L(\lambda)) \\ & = &
\det(\Phi(0,-L,\lambda)Z^-(\lambda),\Phi(0,L,\lambda)Z^+(\lambda))
\rme^{-\mathrm{tr}(A^-(\lambda)P^u_-(\lambda))L} \rme^{\mathrm{tr}(A^+(\lambda)P^s_+(\lambda))L},
\end{eqnarray*}
where the expressions in the last line are all analytic in $\lambda\in\Omega$. Hence, $\tilde{E}_L(\lambda)$ is analytic as claimed.
\end{Proof}

Next, we prove that the system outlined in step~2 in \S\ref{s2.4} has a unique solution.

\begin{Lemma}\label{l2}
There are constants $\eta,C,L_0>0$ so that the following is true for each $L\geq L_0$. Fix $\lambda\in\bar\Omega$, order the eigenvalues of $A^+(\lambda)$ according to their real part, and, starting with $j=1$ and ending at $j=r$, the system
\begin{eqnarray}
V_x = (A(x,\lambda)-\nu^+_j(\lambda)\id)V & \mbox{ on } & 0<x<L \label{e31} \\ \label{e32}
\langle V^+_k(0,\lambda),V(0)\rangle = 0 & \mbox{ for } & k=1,\ldots,j-1 \\ \label{e33}
\langle W^+_j(\lambda),V(L)\rangle = \langle W^+_j(\lambda),V^+_j(\lambda)\rangle && \\ \label{e34}
\langle W^+_k(\lambda),V(L)\rangle = 0 & \mbox{ for } & k=j+1,\ldots,n
\end{eqnarray}
has a unique solution $V(x)=V^+_j(x,\lambda)$ for $x\in[0,L]$.
\end{Lemma}

\begin{Proof}
Firstly, it follows from \citep[Ch.~3 \S8]{CodLev} that there are constants $C,\eta>0$ and $n$ linearly independent solutions $U_j(x,\lambda)$ of
\[
U_x = A(x,\lambda) U, \quad x\geq0
\]
that depend continuously on $\lambda\in\bar\Omega$ such that
\begin{equation}\label{e41}
|U_k(x,\lambda) \rme^{-\nu^+_k(\lambda)x} - V^+_k(\lambda)| \leq C \rme^{-\eta x}, \quad x\geq0
\end{equation}
uniformly in $\lambda\in\bar\Omega$ for $k=1,\ldots,n$. In the remainder of this section, we will use $C$ to denote constants that do not depend on $L$ and $\lambda$.

Next, fix $\lambda\in\bar\Omega$ and order the eigenvalues $\nu^+_k(\lambda)$ of $A^+(\lambda)$ by increasing real part. We will prove the claim by induction over $j$, starting with $j=1$. We claim that we can solve the system (\ref{e31})-(\ref{e34}) iteratively and that the $\ell$th solution is of the form
\begin{equation}\label{e42}
V^+_\ell(0,\lambda) = U_\ell(0,\lambda) + \sum_{k=1}^{\ell-1} a_{k\ell} U_k(0,\lambda) + \rmO(\rme^{-\eta L}).
\end{equation}

For $j=1$, we need to solve
\begin{eqnarray*}
V_x & = & (A(x,\lambda)-\nu^+_1(\lambda)\id)V, \qquad 0\le x\leq L \\
\langle W^+_1(\lambda),V(L)\rangle & = & \langle W^+_1(\lambda),V^+_1(\lambda)\rangle \\
\langle W^+_k(\lambda),V(L)\rangle & = & 0, \qquad\qquad\qquad\qquad\qquad \mbox{ for } k=2,\ldots,n
\end{eqnarray*}
and show that
\[
V^+_1(0,\lambda) = U_1(0,\lambda) + \rmO(\rme^{-\eta L}).
\]
The general solution to the ODE is given by
\[
V(x) = a^c U_1(x,\lambda) \rme^{-\nu^+_1(\lambda)x} + \sum_{k=2}^{n} a^u_k U_k(x,\lambda)
\rme^{-\nu^+_1(\lambda)x} \rme^{(\nu^+_1(\lambda)-\nu^+_k(\lambda))L}
\]
with $a^c\in\C$ and $a^u:=(a^u_2,\ldots,a^u_n)\in\C^{n-1}$ arbitrary. Note that this solution is bounded by $C(|a^c|+|a^u|)$ uniformly in $x\in[0,L]$, independently on $L$, due to the ordering of the eigenvalues $\nu^+_k(\lambda)$. Evaluating at $x=L$, we arrive at
\[
V(L) = a^c (V^+_1(\lambda)+\rmO(\rme^{-\eta L})) + \sum_{k=2}^{n} a^u_k (V^+_k(\lambda)+\rmO(\rme^{-\eta L})).
\]
We can now solve the boundary conditions to find that $a^c=1+\rmO(\rme^{-\eta L})$ and $a^u=\rmO(\rme^{-\eta L})$. In particular, (\ref{e42}) holds for $\ell=1$.

Finally, assume that we solved the system for $\ell=1,\ldots,j-1$ and that (\ref{e42}) holds for all such $\ell$. We consider 
\begin{eqnarray}
V_x & = & (A(x,\lambda)-\nu^+_j(\lambda)\id)V, \qquad 0\leq x\leq L \nonumber \\ \label{p1}
\langle V^+_\ell(0,\lambda),V(0)\rangle & = & 0 \qquad\qquad\qquad\qquad\qquad \mbox{ for } \ell=1,\ldots,j-1 \\ \label{p2}
\langle W^+_j(\lambda),V(L)\rangle & = & \langle W^+_j(\lambda),V^+_j(\lambda)\rangle \\ \label{p3}
\langle W^+_k(\lambda),V(L)\rangle & = & 0, \qquad\qquad\qquad\qquad\qquad \mbox{ for } k=j+1,\ldots,n
\end{eqnarray}
and note that the general solution to the ODE is given by the expression
\[
V(x) = a^c U_j(x,\lambda) \rme^{-\nu^+_j(\lambda)x} + \sum_{k=1}^{j-1} a^s_k U_k(x,\lambda) \rme^{-\nu^+_j(\lambda)x} + \sum_{k=j+1}^{n} a^u_k U_k(x,\lambda) \rme^{-\nu^+_j(\lambda)x} \rme^{(\nu^+_j(\lambda)-\nu^+_k(\lambda))L},
\]
which is again bounded by $C(|a^c|+|a^s|+|a^u|)$ uniformly in $x\in[0,L]$. Evaluating this solution at $x=L$, we obtain
\[
V(L) = a^c (V^+_j(\lambda)+\rmO(\rme^{-\eta L})) + \sum_{k=1}^{j-1} a^s_k (V^+_k(\lambda)+\rmO(\rme^{-\eta L})) \rme^{(\nu^+_k(\lambda)-\nu^+_j(\lambda))L} + \sum_{k=j+1}^{n} a^u_k (V^+_k(\lambda)+\rmO(\rme^{-\eta L})).
\]
Substituting this expression into (\ref{p3}), we can solve those for $a^u=\rmO(\rme^{-\eta L})(a^c,a^s)$. Using the resulting expression for $V(L)$ in (\ref{p2}), we can solve for $a^c\in\C$ and obtain $a^c=1+\rmO(\rme^{-\eta L})a^s$. We therefore find that
\begin{eqnarray*}
V(0) & = & a^c U_j(0,\lambda) + \sum_{k=1}^{j-1} a^s_k U_k(0,\lambda) + \sum_{k=j+1}^{n} a^u_k U_k(0,\lambda) \rme^{(\nu^+_j(\lambda)-\nu^+_k(\lambda))L} \\ & = &
U_j(0,\lambda) + \sum_{k=1}^{j-1} a^s_k U_k(0,\lambda) + \rmO(\rme^{-\eta L}) a^s.
\end{eqnarray*}
Substituting this expression into the remaining boundary conditions (\ref{p1}), we obtain
\begin{equation}\label{p0}
0 = \langle V^+_\ell(0,\lambda),V(0)\rangle = \langle V^+_\ell(0,\lambda), U_j(0,\lambda) \rangle
+ \sum_{k=1}^{j-1} \langle V^+_\ell(0,\lambda),U_k(0,\lambda) \rangle a^s_k + \rmO(\rme^{-\eta L})a^s, \qquad
\ell=1,\ldots,j-1.
\end{equation}
Writing $\mathcal{U}^+_{j-1}(\lambda)$ and $\mathcal{V}^+_{j-1}(\lambda)$ for the matrices in $\C^{n\times(j-1)}$ with columns $U^+_\ell(0,\lambda)$ and $V^+_\ell(0,\lambda)$, respectively, for $\ell=1,\ldots,j-1$, it follows from (\ref{e42}) that there is a matrix $S^+_{j-1}(\lambda)\in\C^{(j-1)\times(j-1)}$ of the form
\begin{equation}\label{e83}
S^+_{j-1}(\lambda) = \begin{pmatrix} 1 & * & \cdots & * \\ 0 & 1 & \ddots & \vdots \\
\vdots & \ddots & \ddots & * \\ 0 & \cdots & 0 & 1\end{pmatrix} + \rmO(\rme^{-\eta L})
\end{equation}
that is bounded uniformly in $\lambda$ such that
\[
\mathcal{V}^+_{j-1}(\lambda) = \mathcal{U}^+_{j-1}(\lambda) S^+_{j-1}(\lambda) + \rmO(\rme^{-\eta L}).
\]
In particular, we see that the matrix
\begin{eqnarray*}
(\langle V^+_\ell(0,\lambda),U_k(0,\lambda) \rangle)_{\ell,k=1,\dots,j-1} & = &
\mathcal{V}^+_{j-1}(\lambda)^* \mathcal{U}^+_{j-1}(\lambda) \\ & = &
\left[\mathcal{U}^+_{j-1}(\lambda) S^+_{j-1}(\lambda)\right]^* \mathcal{U}^+_{j-1}(\lambda)
+ \rmO(\rme^{-\eta L}) \\ & = &
S^+_{j-1}(\lambda)^* \mathcal{U}^+_{j-1}(\lambda)^* \mathcal{U}^+_{j-1}(\lambda) + \rmO(\rme^{-\eta L})
\end{eqnarray*}
is invertible with inverse bounded uniformly in $L$ and $\lambda$, and we can therefore solve (\ref{p0}) uniquely for $a^s\in\C^{j-1}$. Checking the resulting vector $V^+_j(0,\lambda)$, we see that it satisfies (\ref{e42}) with $\ell=j$ as claimed.
\end{Proof}

\begin{Proof}[ of convergence (\ref{et1}) and (\ref{et2}).]
We record from the proof of the preceding lemma that
\begin{eqnarray}
V^+_j(0,\lambda) & = & U^+_j(0,\lambda) + \sum_{k=1}^{j-1} a_{kj} U_k(0,\lambda) + \rmO(\rme^{-\eta L})
\label{e51} \\ \label{e52}
V^+_j(L,\lambda) & = & V^+_j(\lambda) + \sum_{k=1}^{j-1} a_{kj} V^+_k(\lambda) \rme^{(\nu^+_k(\lambda)-\nu^+_j(\lambda))L} + \rmO(\rme^{-\eta L}),
\end{eqnarray}
where the coefficients $a_{kj}$ are bounded uniformly in $L$ and $\lambda$.

Recall that the columns of the matrix $V^+(0,\lambda)$ from \S\ref{s2.2} are given by the vectors $U^+_j(0,\lambda)$. Inspecting (\ref{e51}), we see that there is a uniformly bounded matrix $S^+_L(\lambda)\in\C^{r\times r}$ of the form (\ref{e83}) such that 
\[
V^+(0,\lambda) S^+_L(\lambda) = V^+_L(0,\lambda) + \rmO(\rme^{-\eta L}).
\]
Proceeding analogously for $V^-_L(0,\lambda)$, we see that
\begin{eqnarray*}
\det(V^-(0,\lambda),V^+(0,\lambda)) & = & 
\det(V^-_L(0,\lambda) S^-_L(\lambda)^{-1},V^+_L(0,\lambda)S^+_L(\lambda)^{-1}) + \rmO(\rme^{-\eta L}) \\ & = &
\det(V^-_L(0,\lambda),V^+_L(0,\lambda)) (1+\rmO(\rme^{-\eta L})) + \rmO(\rme^{-\eta L})
\end{eqnarray*}
and therefore
\[
|\mathcal{E}(\lambda) - \tilde{\mathcal{E}}_L(\lambda)|
= |\det(V^-(0,\lambda),V^+(0,\lambda)) - \det(V^-_L(0,\lambda),V^+_L(0,\lambda))|
\leq C \rme^{-\eta L}
\]
as claimed.

Next, we recall that the matrices $C^+_L(\lambda)$ and $\tilde{C}^+_L(\lambda)$ are solutions to the systems (\ref{e12}) and (\ref{e81}), respectively, which we write as
\begin{equation}\label{e82}
V^+(L,\lambda) C^+_L(\lambda) = Z^+(\lambda) + \rmO(\rme^{-\eta L}), \qquad
V^+_L(L,\lambda) \tilde{C}^+_L(\lambda) = Z^+(\lambda).
\end{equation}
Equation (\ref{e52}) implies that there is a uniformly bounded matrix $T^+_L(\lambda)\in\C^{r\times r}$ of the form (\ref{e83}) such that
\[
V^+(L,\lambda) T^+_L(\lambda) = V^+_L(L,\lambda) + \rmO(\rme^{-\eta L}).
\]
Hence, we can rewrite the second equation in (\ref{e82}) as
\[
(V^+(L,\lambda) T^+_L(\lambda) + \rmO(\rme^{-\eta L})) \tilde{C}^+_L(\lambda) = Z^+(\lambda) 
\]
and, using the first equation in (\ref{e82}), we arrive at
\[
(V^+(L,\lambda) T^+_L(\lambda) + \rmO(\rme^{-\eta L})) \tilde{C}^+_L(\lambda) = 
V^+(L,\lambda) C^+_L(\lambda) + \rmO(\rme^{-\eta L}).
\]
Using the exponential convergence of $V^+(L,\lambda)$ to the full-rank matrix $V^+_\infty(\lambda)$, whose columns are formed by the eigenvectors $V^+_j(\lambda)$ of $A^+(\lambda)$, we conclude that
\[
C^+_L(\lambda) = T^+_L(\lambda) \tilde{C}^+_L(\lambda) + \rmO(\rme^{-\eta L})
\]
and therefore
\[
\det C^+_L(\lambda) = \det \tilde{C}^+_L(\lambda) + \rmO(\rme^{-\eta L}).
\]
This implies $|\mathcal{C}_L(\lambda)-\tilde{\mathcal{C}}_L(\lambda)|\leq C\rme^{-\eta L}$ and completes the proof of (\ref{et1}). The remaining estimate (\ref{et2}) now follows easily.
\end{Proof}

\begin{Remark}\label{r:1}
We emphasize that our proof demonstrates that the algorithm proposed in \S\ref{s2.4} is well-posed with bounds on the solutions that are independent of $L$ and $\lambda$ even when we compute solutions for $j=1,\ldots,k$ for $k>r$ as long as the eigenvalues $\nu^+_j(\lambda)$ can be ordered according to increasing real part. Indeed, the expressions for the solutions given above depend only on $\nu^+_j(\lambda)-\nu^+_k(\lambda)$ with $j<k$ so that solutions are always bounded. This observation is what makes our algorithm suitable for computing Evans functions  that are extended analytically across the essential spectrum: see \S\ref{s:kdv} for an example of such a computation.
\end{Remark}

% --------------------------------------------------------------------------------------

\section{Implementation and benchmark computations}\label{s3}

In this section, we show how our algorithm performs on various test systems. In particular, we compare the computational costs of the algorithm proposed here with continuous orthogonalization and demonstrate the accuracy of our algorithm in computing the correct winding number for detecting roots of the Evans function.

% --------------------------------------------------------------------------------------

\subsection{Implementation}\label{s3.1}

In the following sections, we will compare continuous orthogonalization as developed in \citep{HuZ2,Z5} to the algorithm presented in \S\ref{s2.4}. We used the implementation of continuous orthogonalization in the Matlab version of \textsc{stablab} (see \citep{STABLAB}) and solved the resulting nonlinear ODE systems using Matlab's \texttt{ode15s} solver, which is an adaptive implicit ODE solver. We used three different Matlab implementations for the algorithm outlined in \S\ref{s2.4}. The first version uses the Matlab function \texttt{bvp5c}, which is a finite-difference code for solving nonlinear boundary-value problems that implements the four-stage Lobatto IIIa formula and solves the resulting algebraic system directly. The second version uses the Matlab script \texttt{bvp6c} \citep{HM}, which uses a 6th-order interpolant. Thirdly, we also wrote a Matlab script (referred to as \texttt{bvpcheb} below) that solves linear boundary-value problems using differentiation matrices based on collocation with Chebyshev polynomials. In each of these four implementations, we provided the analytic Jacobian of the ODE equations and the boundary conditions to improve computational speed.

% --------------------------------------------------------------------------------------

\subsection{Coupled Nagumo system}\label{s:cn}

\begin{figure}
\centering
\includegraphics[width=0.3\textwidth]{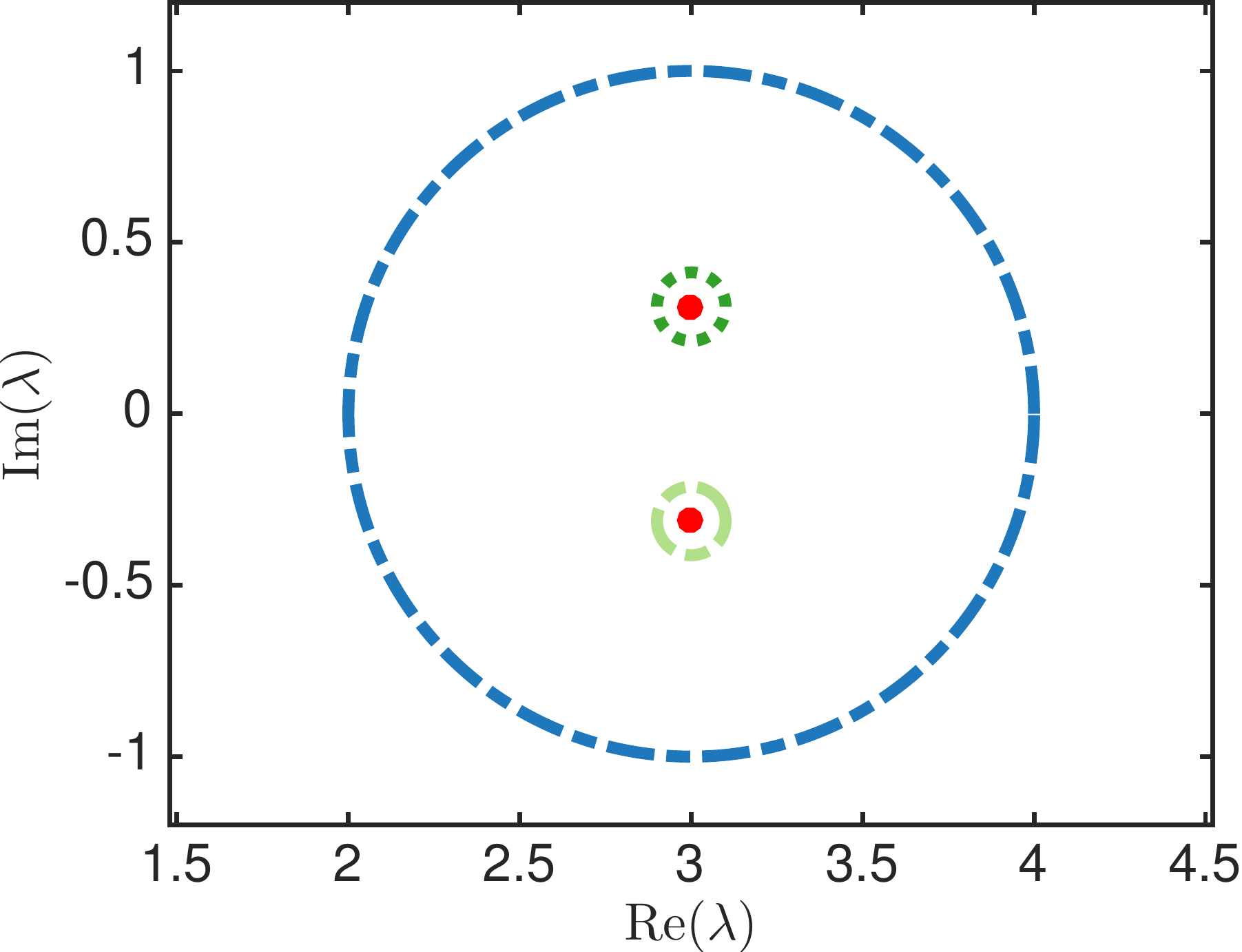} \qquad
\includegraphics[width=0.3\textwidth]{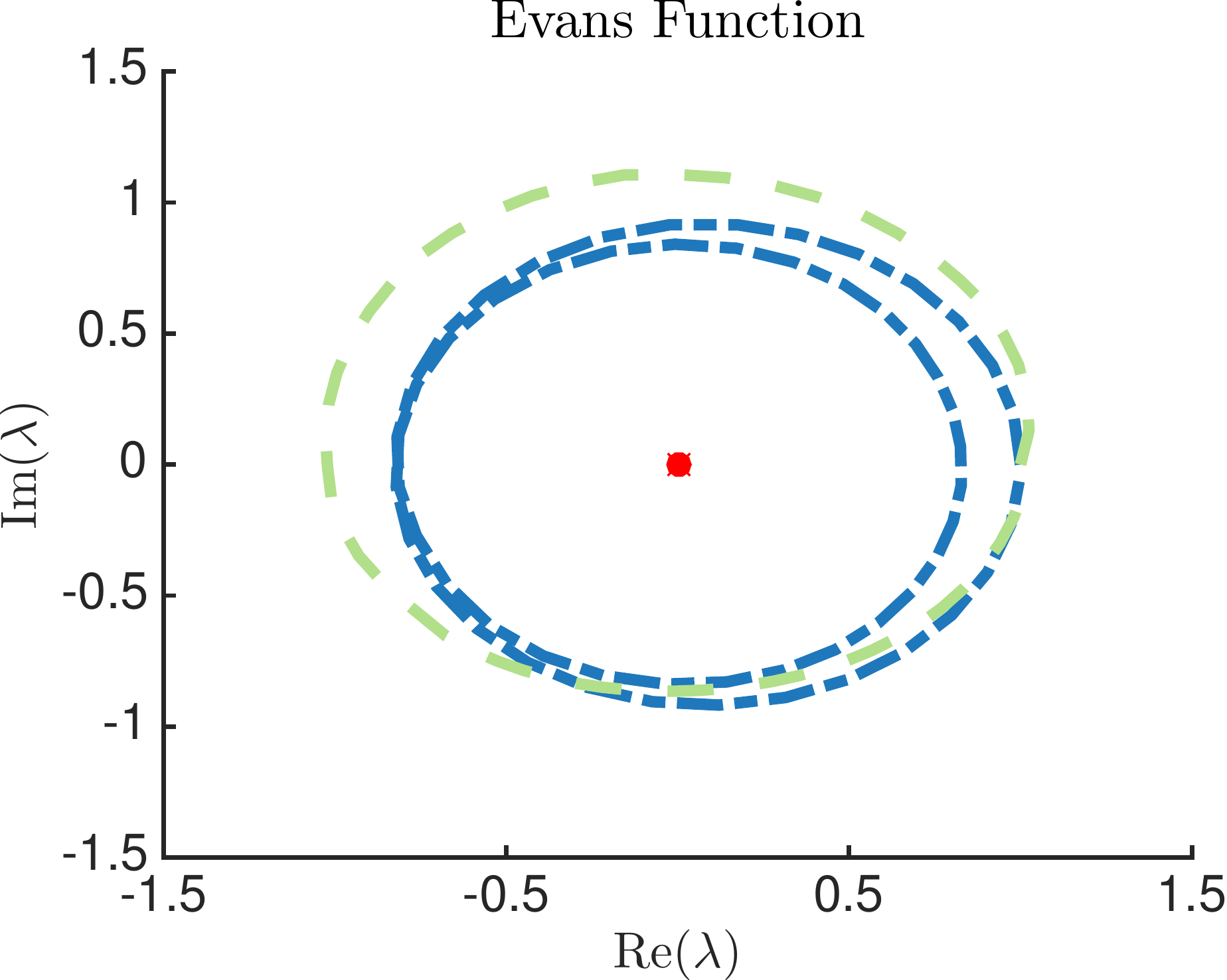} \qquad
\includegraphics[width=0.3\textwidth]{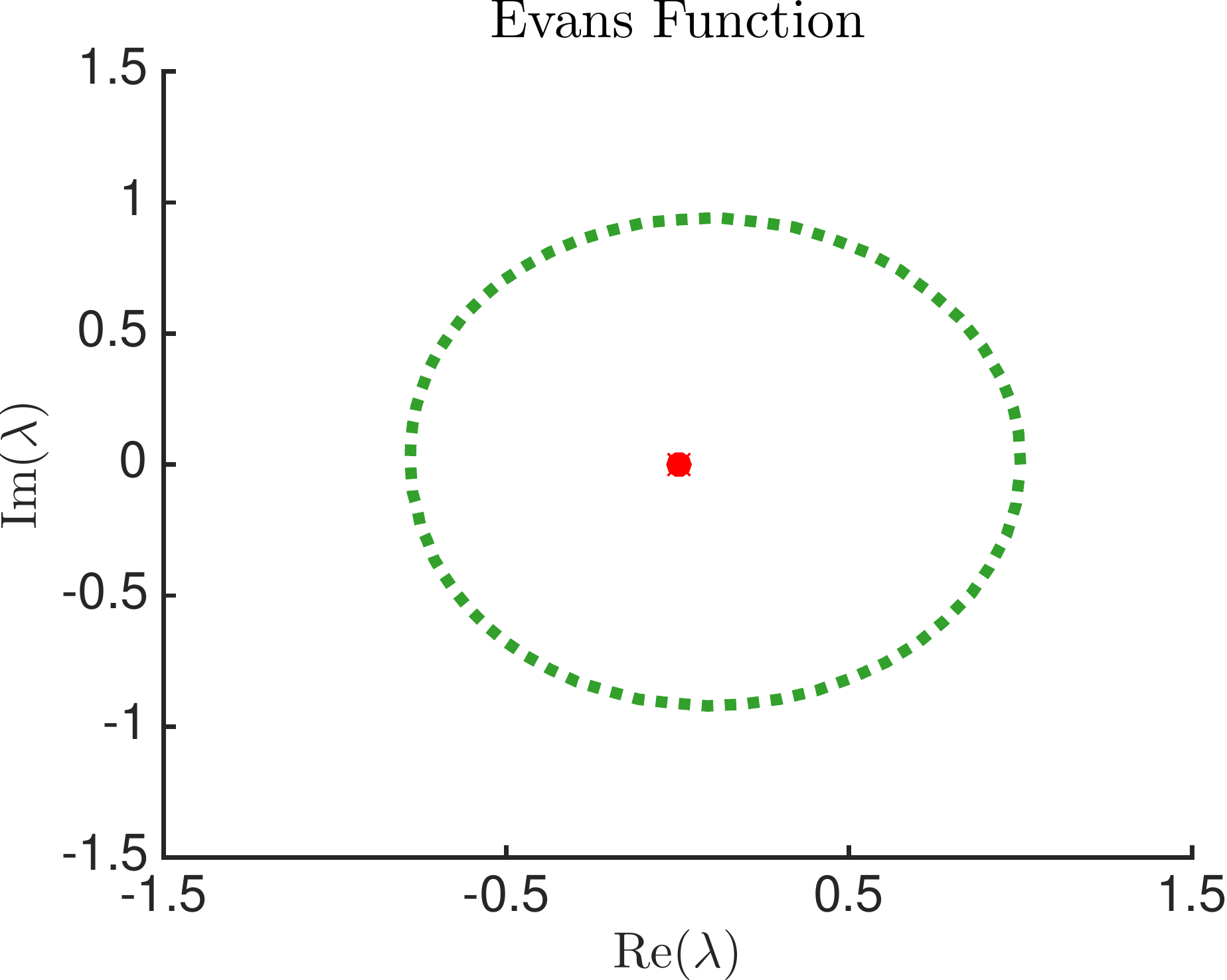} 
\caption{Shown are Evans-function computations for the coupled Nagumo system (\ref{e:node}) with $a=0.1$ and $b=-1$ using our algorithm. The exact eigenvalues are $\lambda=3\pm\rmi/\sqrt{10}$ and $\lambda=\pm\rmi/\sqrt{10}$. The left panel indicates the three contours along which we compute the Evans function, namely a circle of radius 1 centered at 3 (dark blue dash-dotted line) and circles of radius 0.1 centered at $3\pm5/16\rmi$ (dark green dotted line and light green dashed line, respectively). The images of the contours under the Evans function are shown in the center and right panels using the same line styles (the red square indicates the origin).}
\label{f:n1}
\end{figure}

The scalar Nagumo equation takes the form
\[
u_t = u_{xx} - u (1-u^2),
\]
where $x\in\R$. This system admits the localized stationary solution $u_*(x)=\sqrt{2}\sech(x)$. The resulting eigenvalue problem is given by
\[
\lambda u = u_{xx} - (1-3u_*(x)^2) u
\]
which gives the simple eigenvalues $\lambda=3$ and $\lambda=0$ plus the essential spectrum $\{\lambda\in\R:\;\lambda\leq-1\}$. We consider the coupled eigenvalue problem
\begin{eqnarray}
\lambda u & = & u_{xx} - (1-3u_*(x)^2) u + a v \label{e:npde} \\ \nonumber
\lambda v & = & v_{xx} - (1-3u_*(x)^2) v + b u,
\end{eqnarray}
which has eigenvalues at $\lambda=3\pm\sqrt{ab}$ and $\lambda=\pm\sqrt{ab}$ and essential spectrum at $\{\Re\lambda\leq-1, \; \Im\lambda=\pm\sqrt{ab}\}$. Written as a first-order system, the eigenvalue problem (\ref{e:npde}) becomes
\begin{equation}\label{e:node}
U_x = \begin{pmatrix} 0 & 1 & 0 & 0 \\ \lambda+1-3u_*(x)^2 & 0 & -a & 0\\
0 & 0 & 0 & 1 \\ -b & 0 & \lambda+1-3u_*(x)^2 & 0 \end{pmatrix} U
=: A(x,\lambda) U, \qquad U\in\C^4.
\end{equation}
For our computations, we chose $a=0.1$ and $b=-1$ to demonstrate that our algorithm works well in the case where the asymptotic matrix $A_0(\lambda)=\lim_{x\to\pm\infty}A(x,\lambda)$ has complex conjugate eigenvalues. Figure~\ref{f:n1} shows the Evans function computed on a circle of radius $1$ centered at $\lambda=3$ and on circles of radius $0.1$ centered at $\lambda=3\pm\frac{5}{16}\rmi$.

% --------------------------------------------------------------------------------------

\subsection{Korteweg--de Vries equation}\label{s:kdv}

\begin{figure}
\centering 
\includegraphics[width=0.24\textwidth]{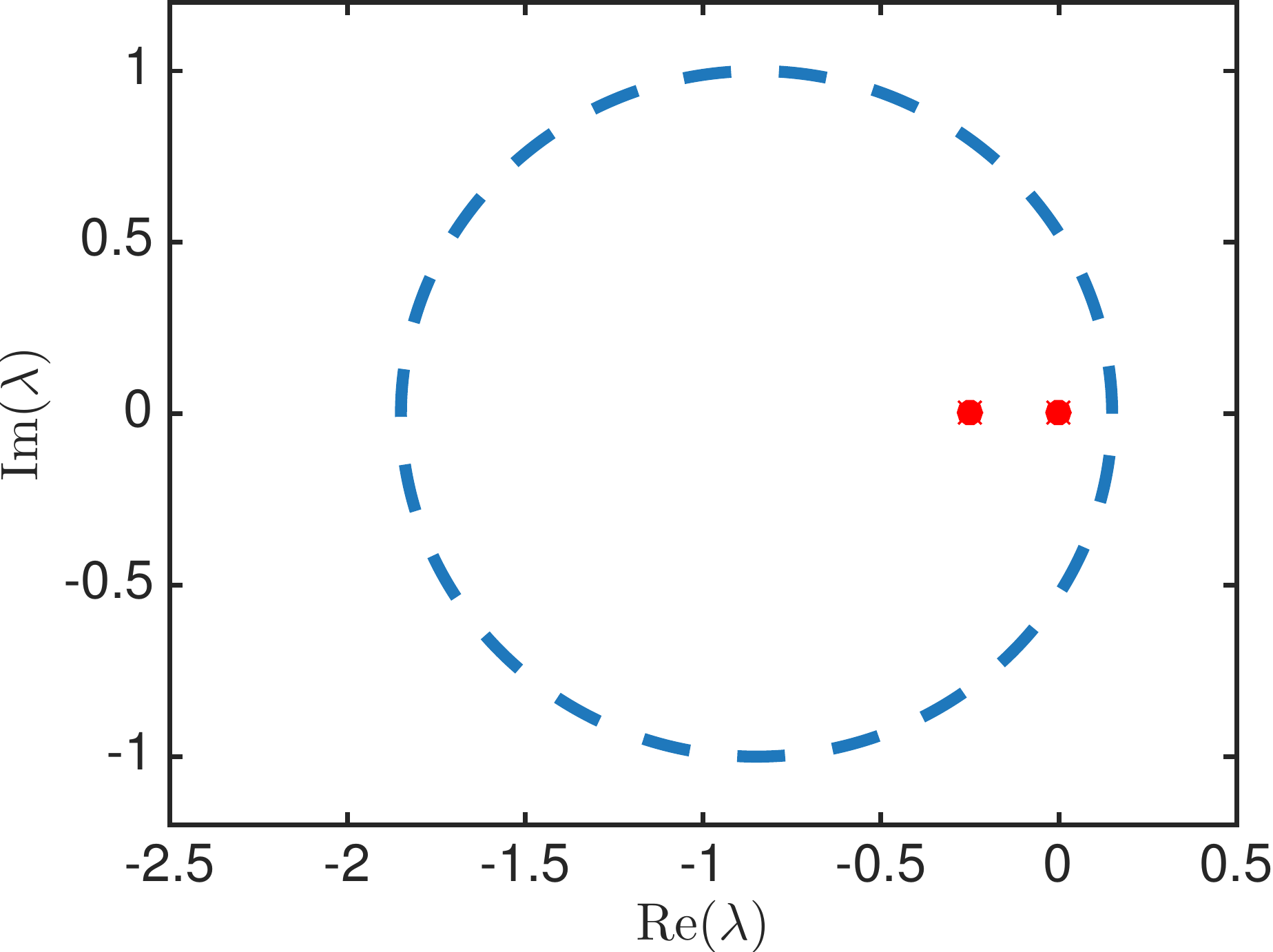} 
\includegraphics[width=0.24\textwidth]{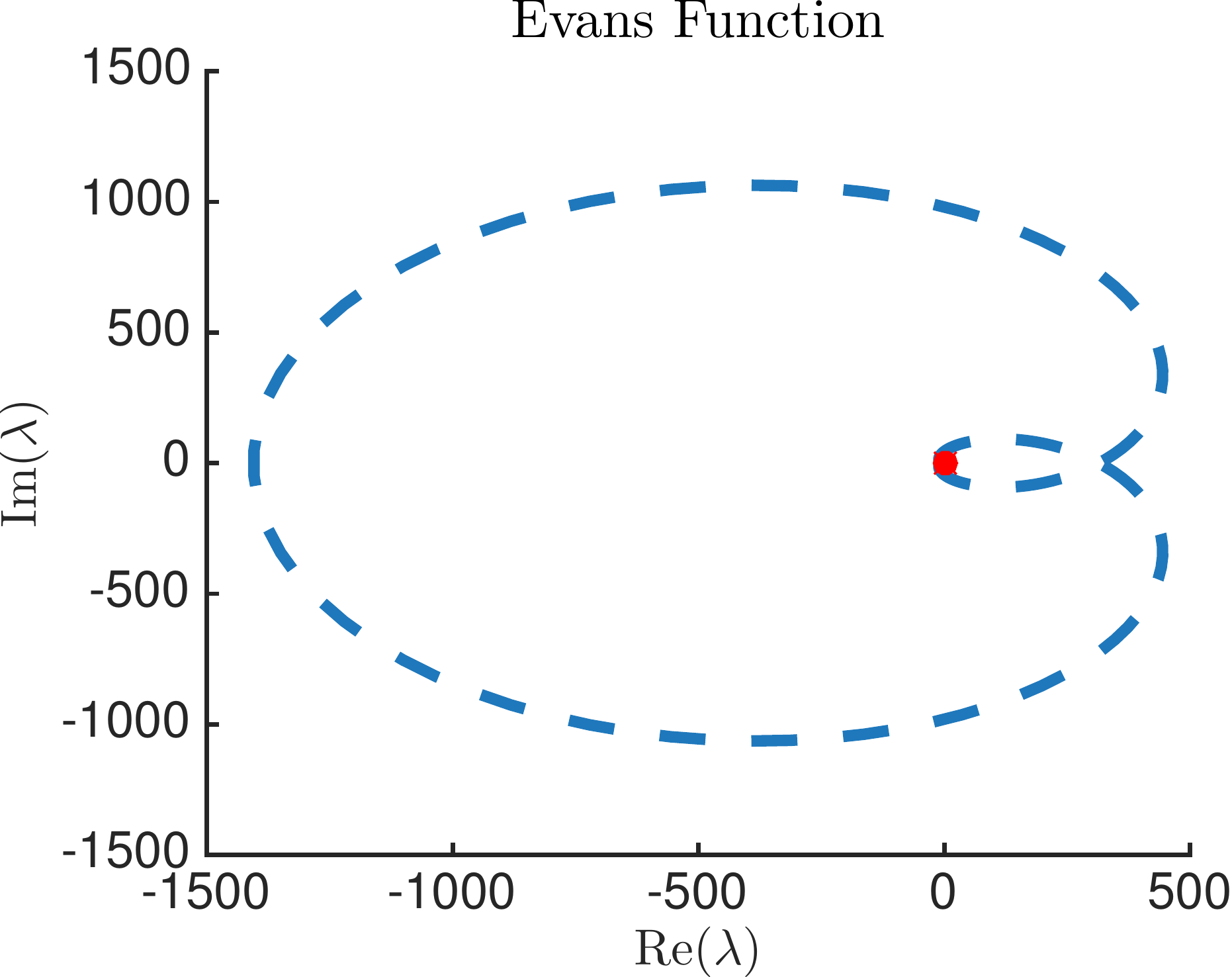} 
\includegraphics[width=0.24\textwidth]{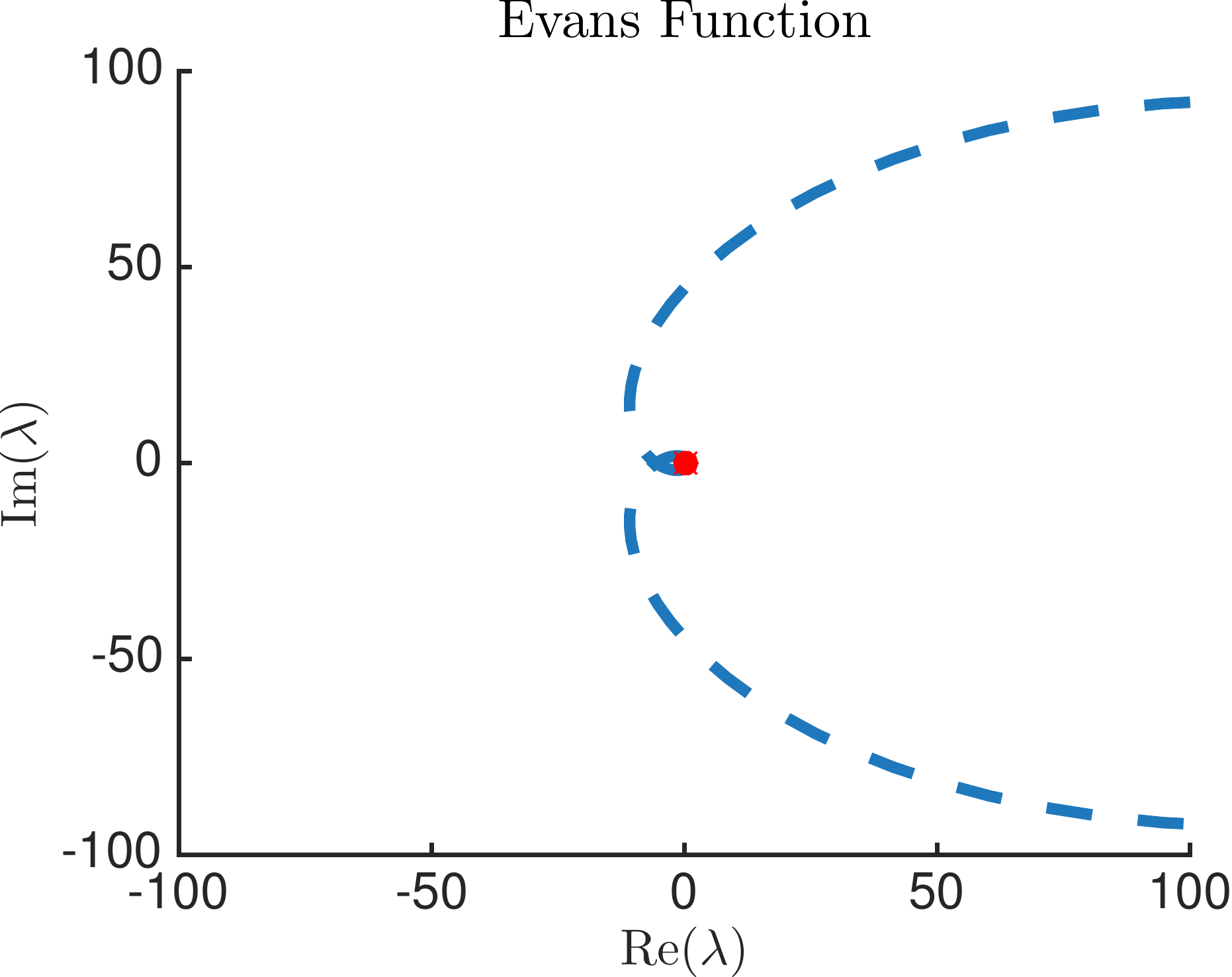} 
\includegraphics[width=0.24\textwidth]{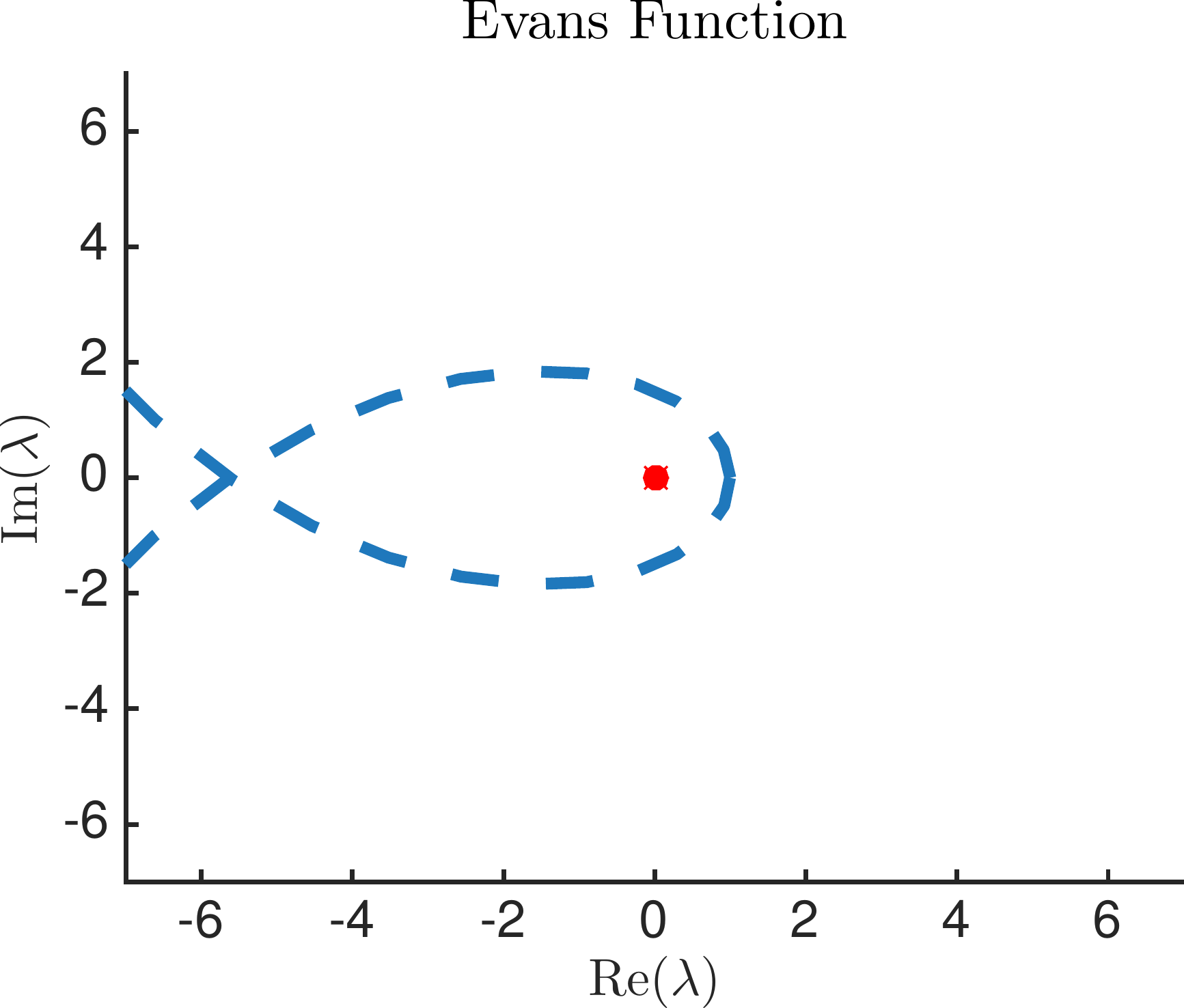} 
\caption{Shown are Evans-function computations for the Korteweg--de Vries equation (\ref{e:kdvode}) with $p=3.95$ and $c=5$. The Evans function has a simple root at $\lambda\approx-0.25$ and a double root at $\lambda=0$. The left panel contains the contour, a circle of radius 1 centered at $\lambda=-0.85$, together with the anticipated roots indicated as red squares. The right three panels show the image of the contour under the Evans function with the origin indicated as a red square: the winding number is 3 as expected.}
\label{f:kdv1}
\end{figure}

\begin{figure}
\centering
\includegraphics[width=0.24\textwidth]{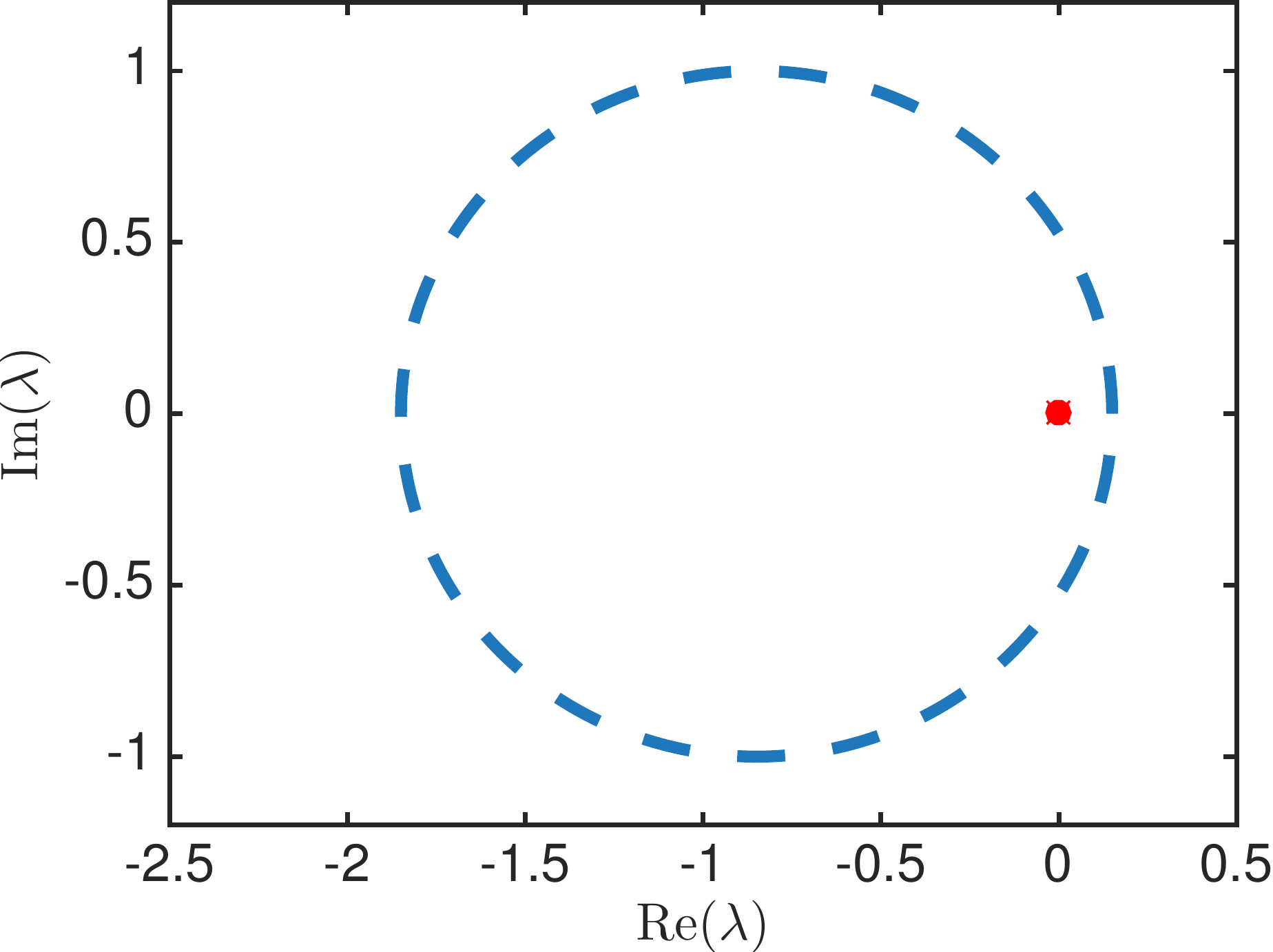} 
\includegraphics[width=0.24\textwidth]{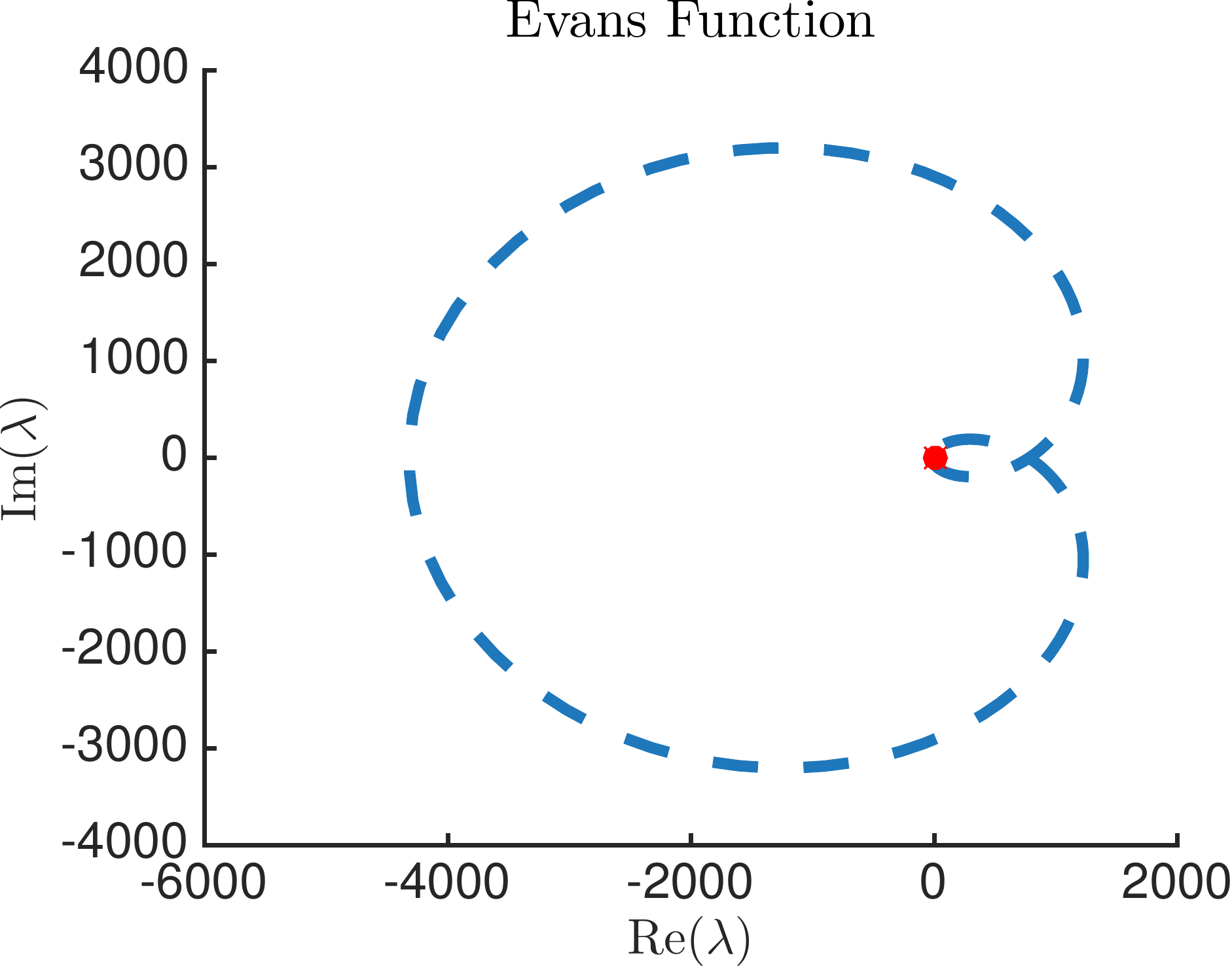} 
\includegraphics[width=0.24\textwidth]{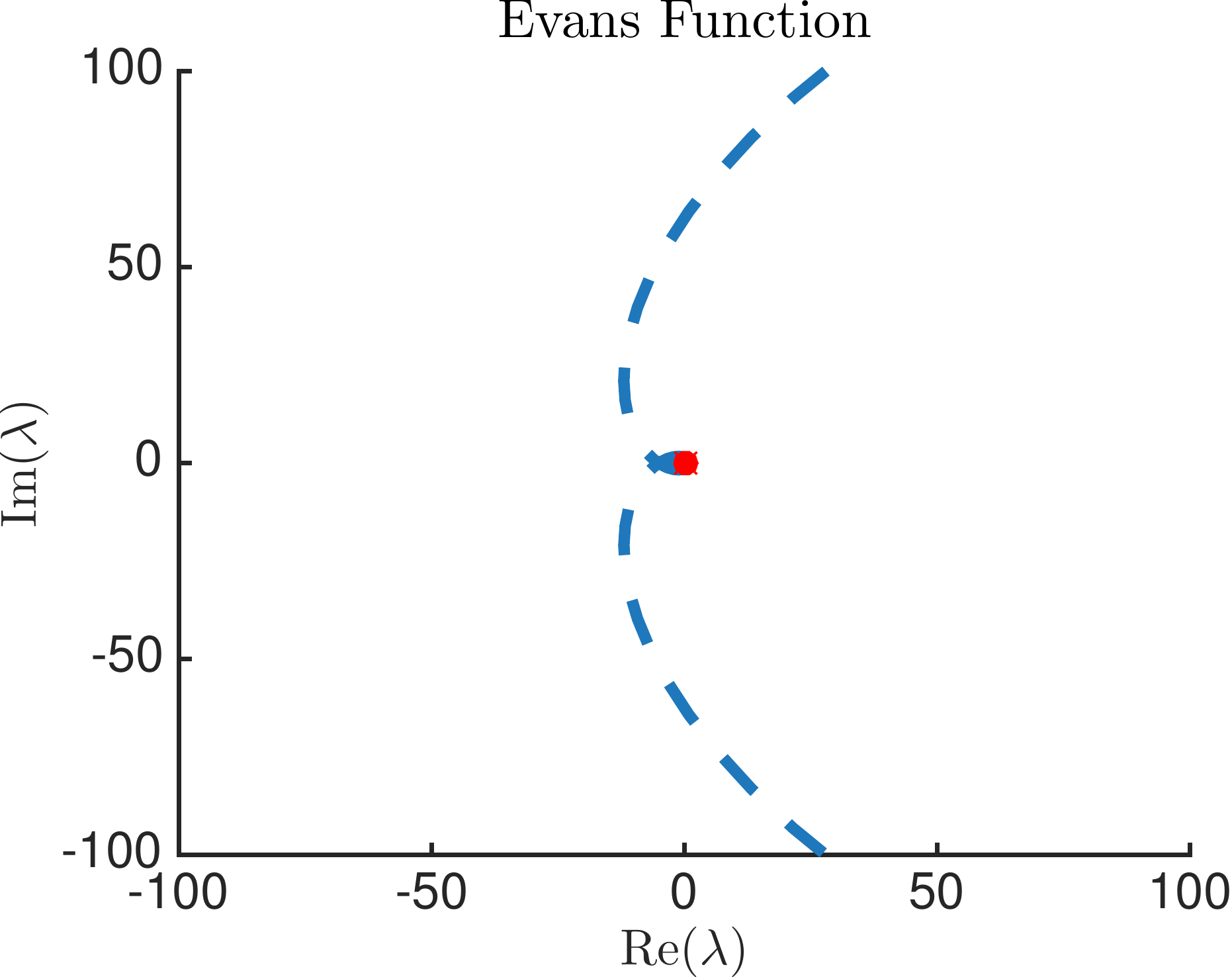} 
\includegraphics[width=0.24\textwidth]{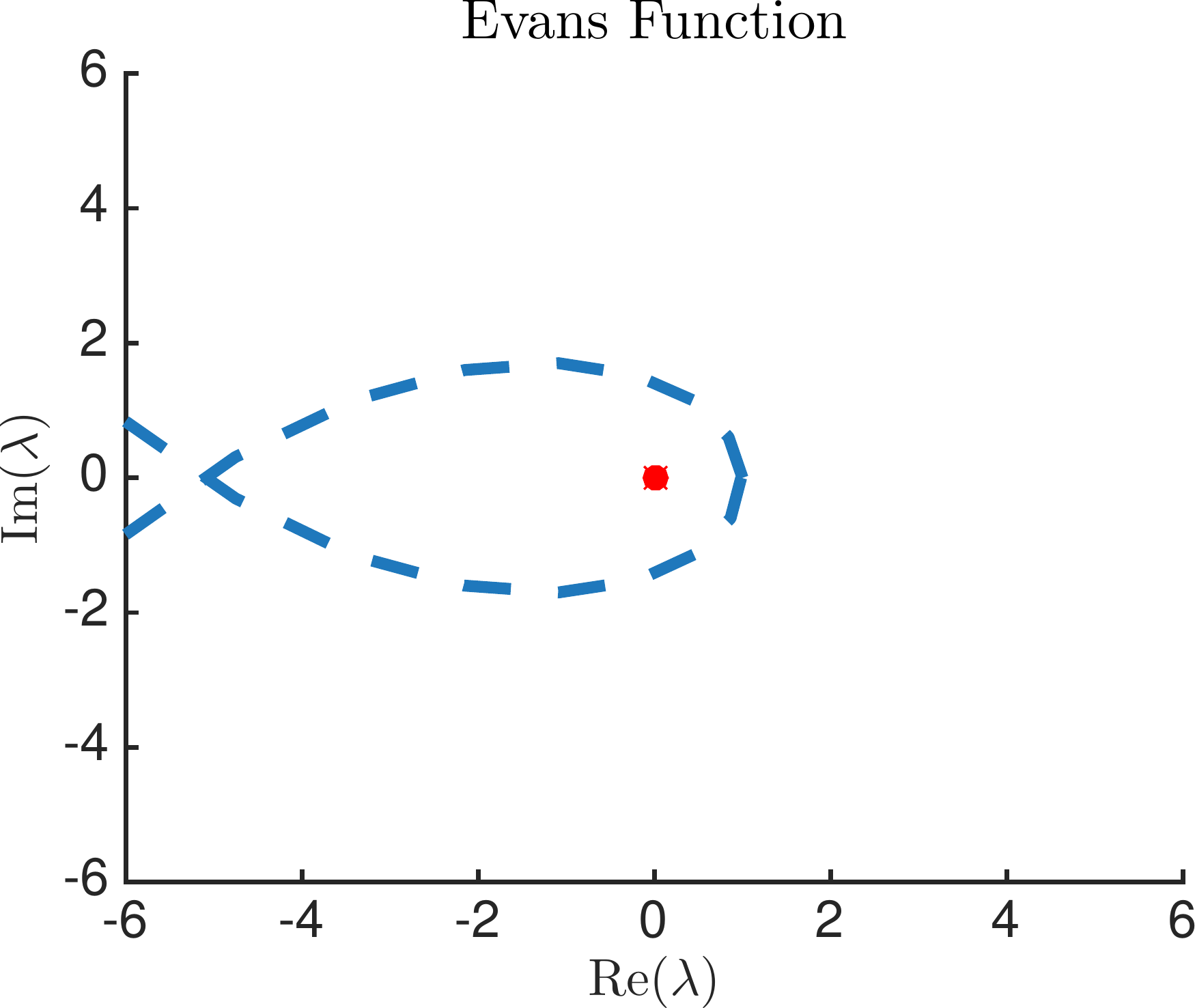} 
\caption{Shown are Evans-function computations for the Korteweg--de Vries equation (\ref{e:kdvode}) with $p=4$ and $c=5$. The Evans function has a triple root at $\lambda=0$. The left panel contains the contour, a circle of radius 1 centered at $\lambda=-0.85$, together with the anticipated root indicated as a red square. The right three panels show the image of the contour under the Evans function with the origin indicated as a red square: the winding number is again 3 as expected.}
\label{f:kdv2}
\end{figure}

\begin{figure}
\centering
\includegraphics[width=0.24\textwidth]{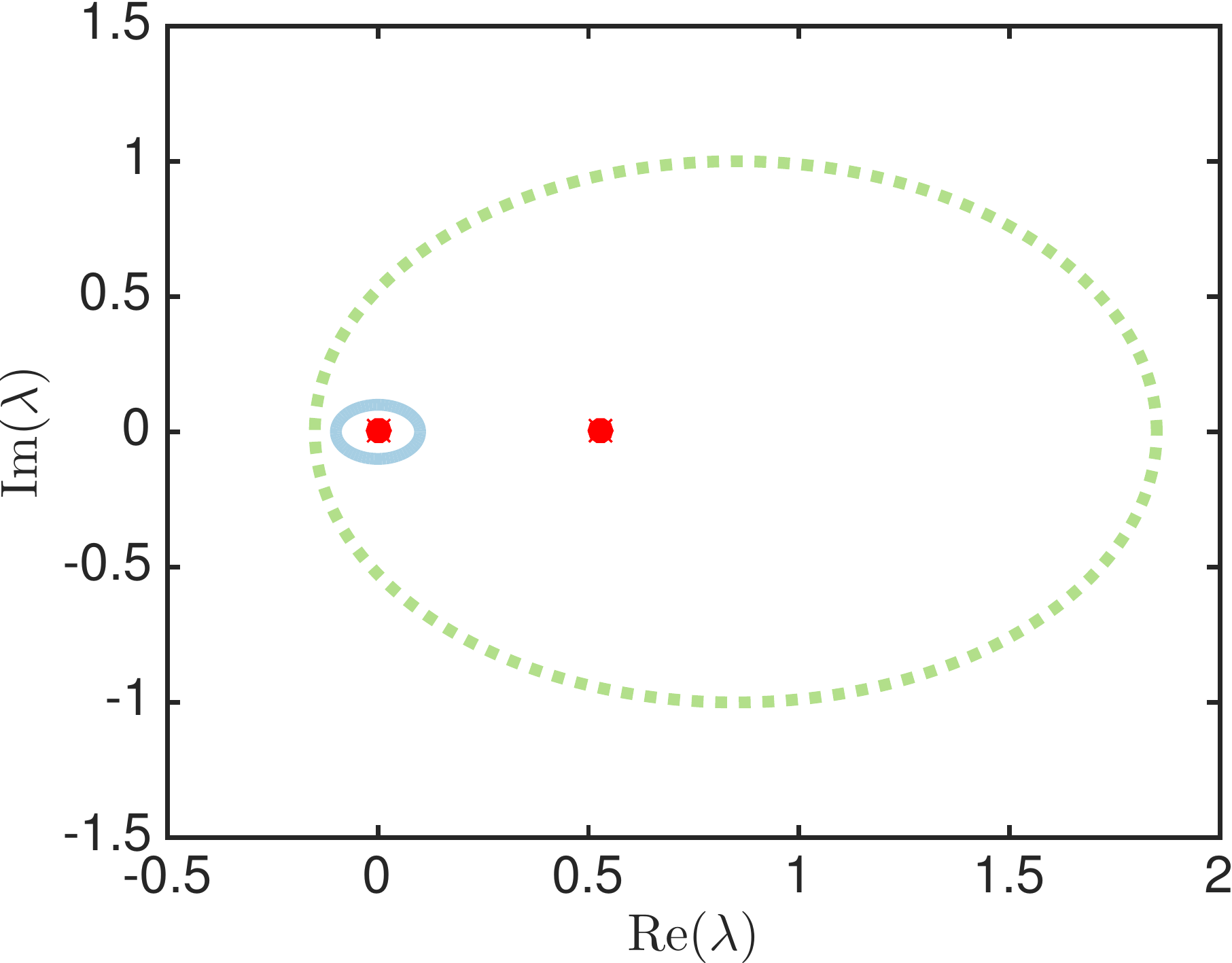} 
\includegraphics[width=0.24\textwidth]{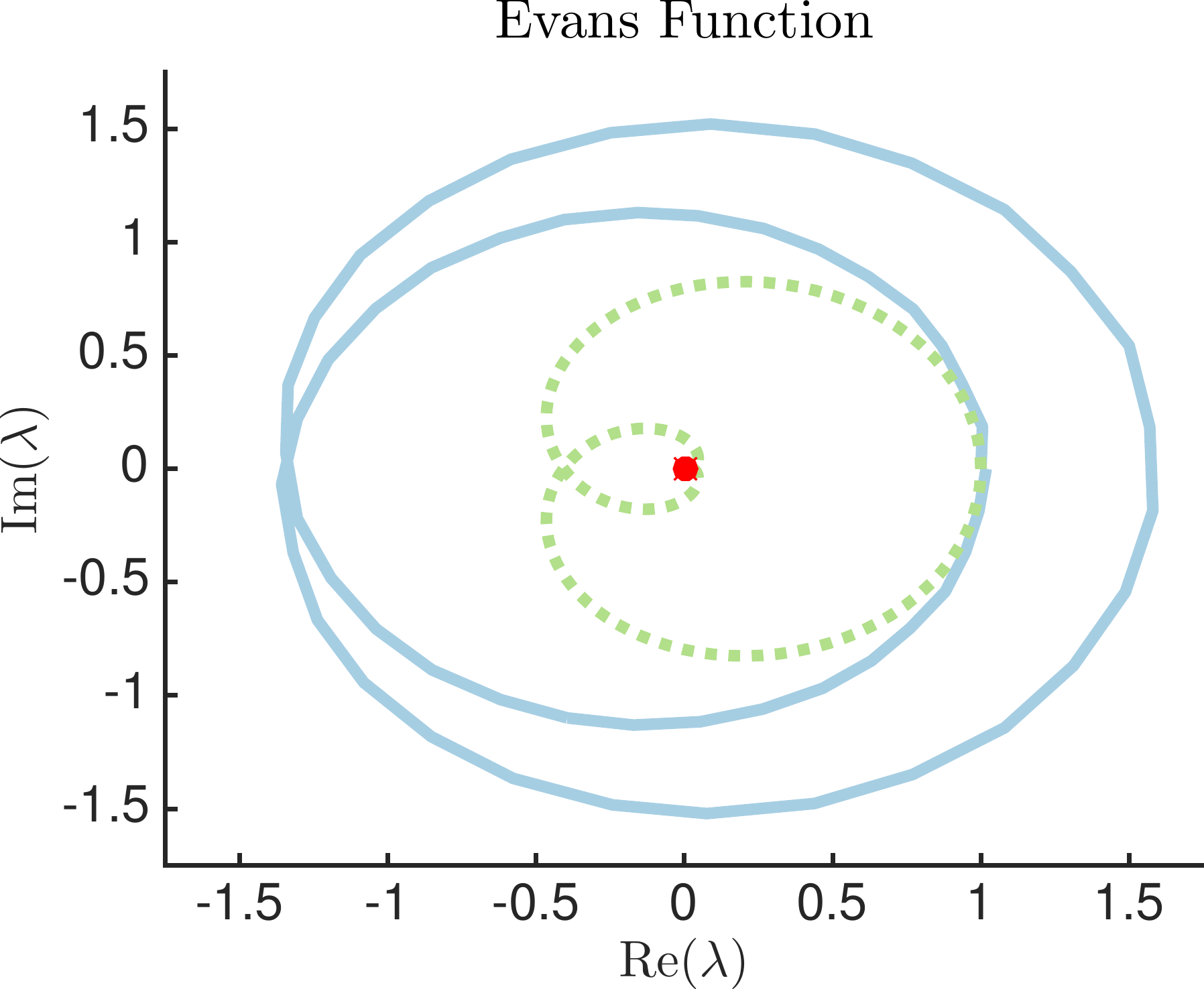} 
\includegraphics[width=0.24\textwidth]{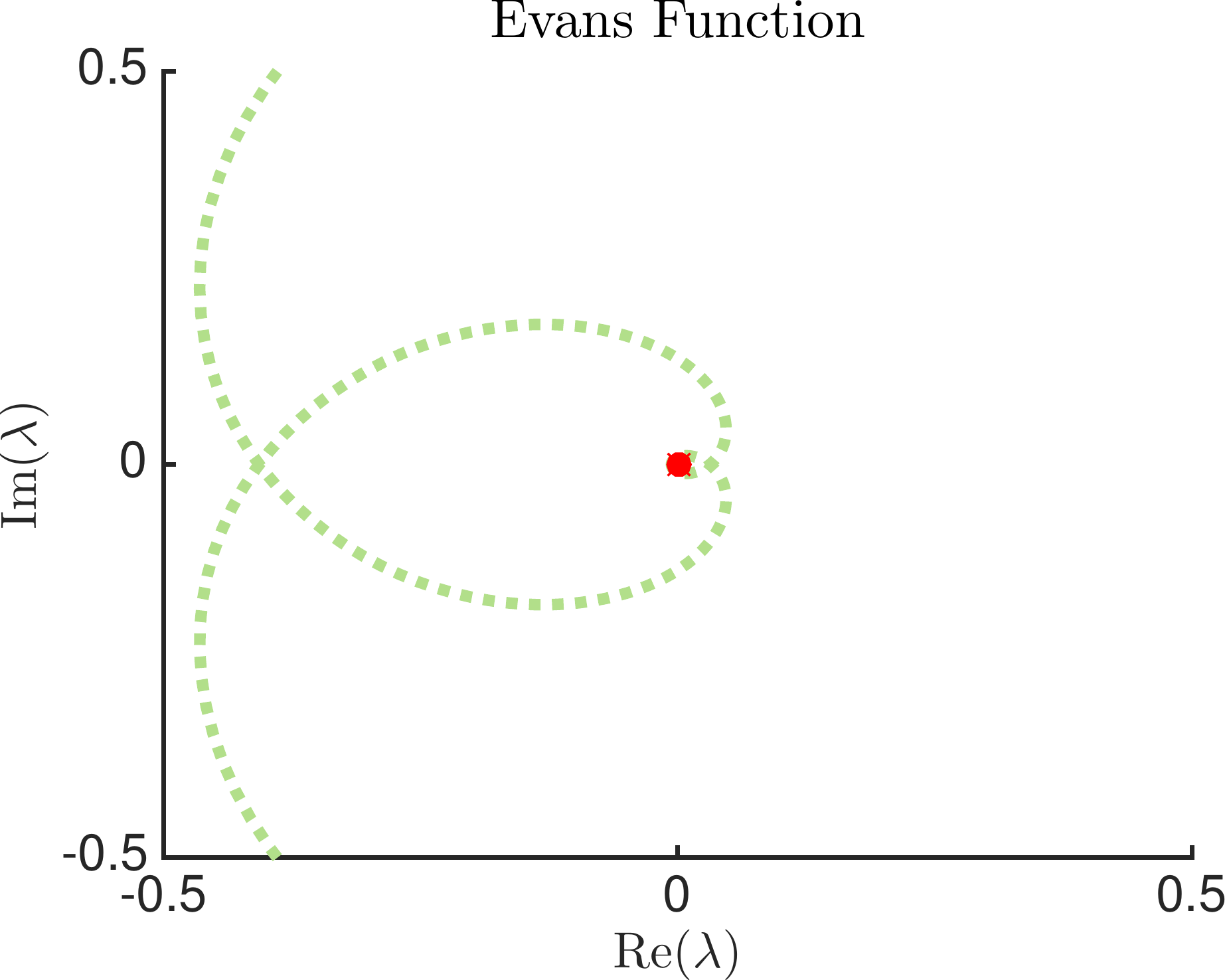} 
\includegraphics[width=0.24\textwidth]{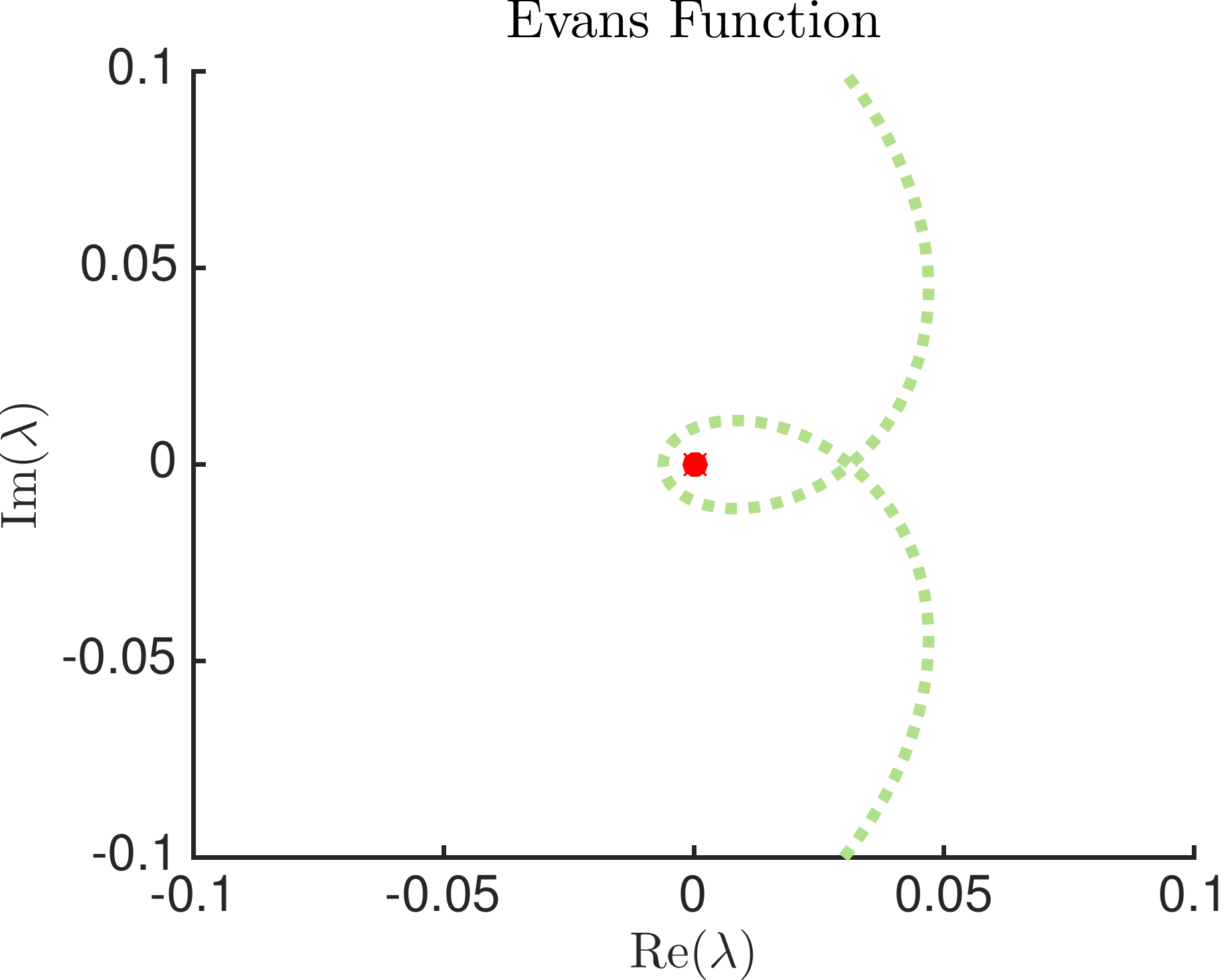} 
\caption{Shown are Evans-function computations for the Korteweg--de Vries equation (\ref{e:kdvode}) with $p=4.1$ and $c=5$. The Evans function has a simple root at $\lambda\approx0.53$ and a double root at $\lambda=0$. The left panel contains the two contours we used, namely a circle of radius 1 centered at $\lambda=0.85$ (dashed light green) and circle of radius 0.1 centered at the origin (solid light blue) together with the anticipated roots indicated as red squares. The right three panels show the image of these contours under the Evans function with the origin indicated as a red square: The winding number of the solid light blue contour is two, while the dashed light green contour gives a winding number of three as expected.}
\label{f:kdv3}
\end{figure}

Next, we consider the Korteweg--de Vries equation
\[
u_t + u_{xxx} + \frac{1}{p+1} (u^{p+1})_x = 0,
\]
where $p>2$. Transforming into the moving coordinate $x\mapsto x-ct$, we arrive at the system
\begin{equation}\label{e:kdv}
u_t + u_{xxx} - cu_x + \frac{1}{p+1} (u^{p+1})_x = 0,
\end{equation}
which admits the stationary solitary waves
\[
u_*(x) = \left[\frac{1}{2}c(p+2)(p+1)\right]^{\frac{1}{p}} \sech\left(\frac{\sqrt{c}p}{2}x\right)^{\frac{2}{p}}.
\]
Linearizing (\ref{e:kdv}) about $u_*(x)$, we obtain the eigenvalue problem
\begin{equation}\label{e:kdepde}
\lambda u + u_{xxx} - cu_x + (u_*(x)^p u)_x = 0,
\end{equation}
which we rewrite as the first-order system
\begin{equation}\label{e:kdvode}
U_x = \begin{pmatrix} 0 & 1 & 0 \\ 0 & 0 & 1 \\ -\lambda - pu_*(x)^{p-1}u_*^\prime(x) & c-u_*(x)^p & 0 \end{pmatrix}
=: A(x,\lambda) U, \qquad U\in\C^3.
\end{equation}
Note that the matrix $A_0(\lambda):=\lim_{x\to\pm\infty}A(x,\lambda)$ fails to be hyperbolic for $\lambda\in\rmi\R$ as the imaginary axis consists of essential spectrum for the PDE linearization. As shown in \citep{BDG,PW}, the Evans function can be extended analytically across the imaginary axis, and the resulting function has a double root at the origin independently of $c$ and $p$, and a simple root on the real axis that crosses from the left into the right half-plane as the parameter $p$ crosses through $p=4$. In Figures~\ref{f:kdv1}-\ref{f:kdv3}, we demonstrate that our numerical Evans-function algorithm can be set up to correctly compute the extended Evans function. In particular, our computations indicate that our algorithm picks up the movement of this root as well as the location of the double root at $\lambda=0$ correctly for $p=3.95,4,4.1$.

% --------------------------------------------------------------------------------------

\subsection{Swift--Hohenberg equation}\label{s:sh}

\begin{figure}
\centering
\includegraphics[height=0.2\textwidth]{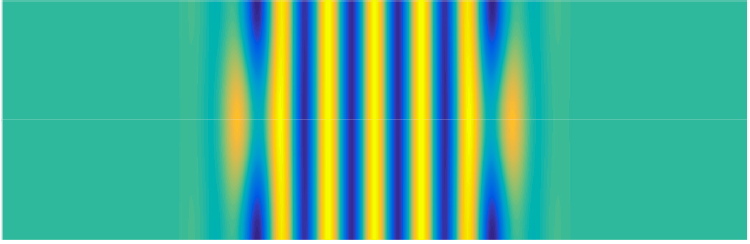}
\caption{Shown is a stationary solution $u_*(x,y)$ of the Swift--Hohenberg equation (\ref{e:sh}) with parameters $\mu=0.675$ and $\nu=2$ on the domain $[-50,50]\times[0,2\pi]$ with periodic boundary conditions in the transverse $y$-direction. Note that this solution is localized in the $x$-direction.}
\label{f:sh1}
\end{figure}

The third example we consider is the planar Swift--Hohenberg equation
\begin{equation}\label{e:sh}
u_t = -(1+\Delta)^2u - \mu u + \nu u^3 - u^5
\end{equation}
for $(x,y)\in[-50,50]\times[0,2\pi]$ with periodic boundary conditions in the $y$-direction. This equation admits the stationary solution $u_*(x,y)$ shown in Figure~\ref{f:sh1} (see \citep{Avitabile2010}).

\begin{figure}
\centering
\includegraphics[width=0.33\textwidth]{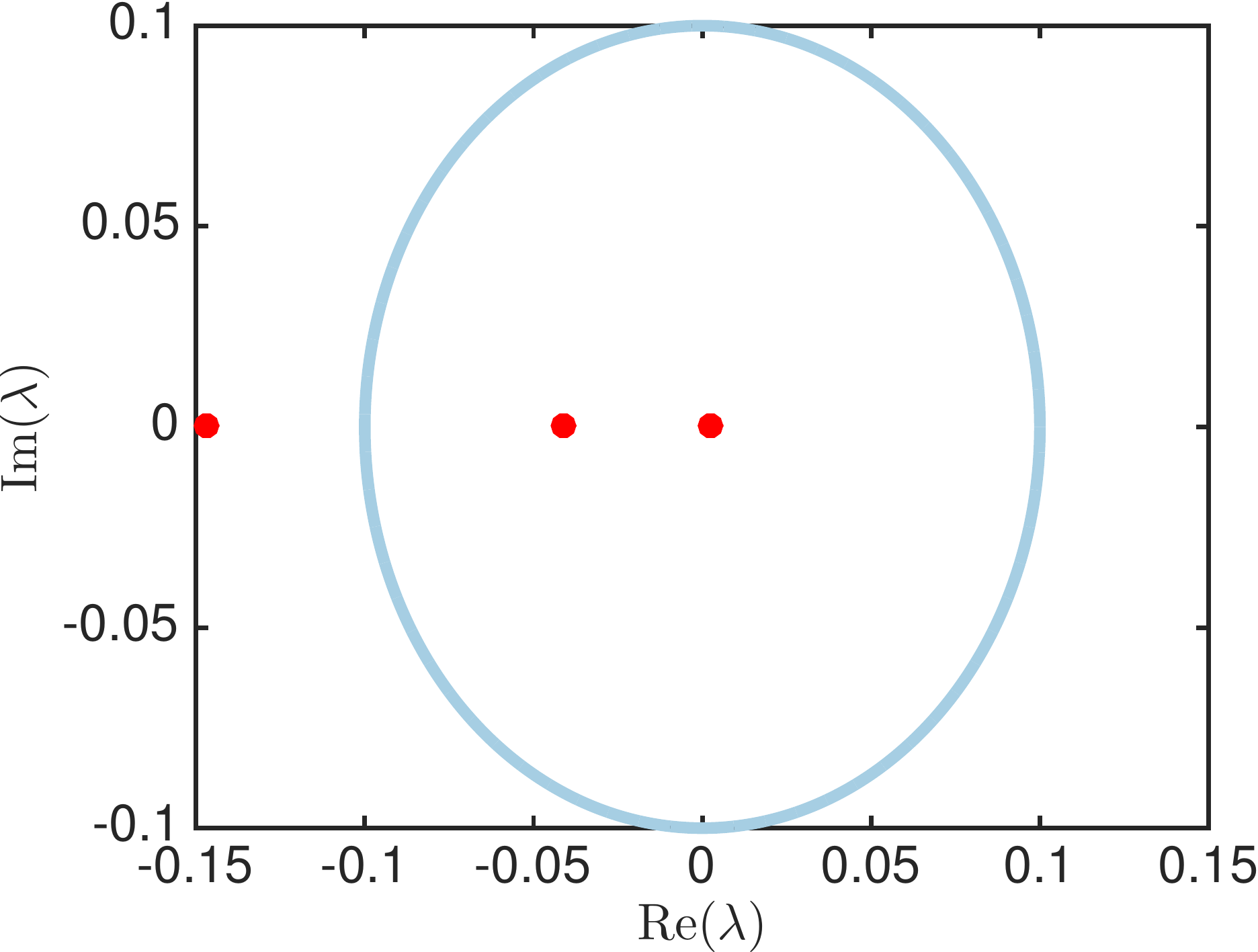} \qquad
\includegraphics[width=0.3\textwidth]{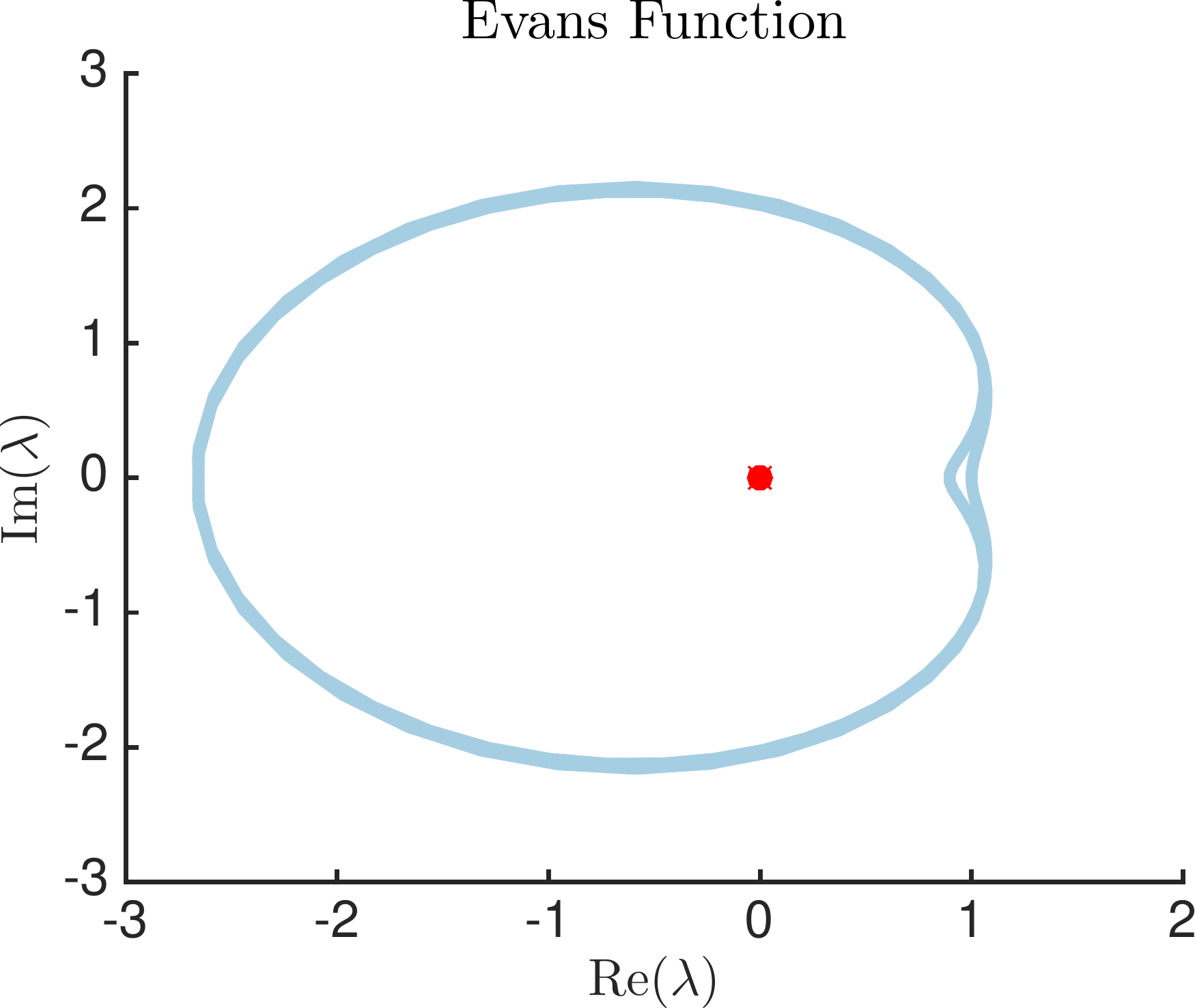} 
\caption{Shown is the Evans function of the Swift--Hohenberg equation (\ref{e:shode}) for $\mu=0.675$ and $\nu=2$. The left panel contains the contour, a circle of radius 1 centered at the origin that surrounds precisely the eigenvalues $\lambda=0,-0.0412$ (red squares); also shown is the eigenvalue $\lambda=-0.1472$ (red square) outside the circle. The right panel shows the image of contour under the Evans function with the original marked as a red square: the winding number is 2, reflecting the combined multiplicity of the eigenvalues enclosed by the circle in the left panel.}
\label{f:sh2}
\end{figure}

\begin{figure}
\centering 
\includegraphics[width=0.19\textwidth]{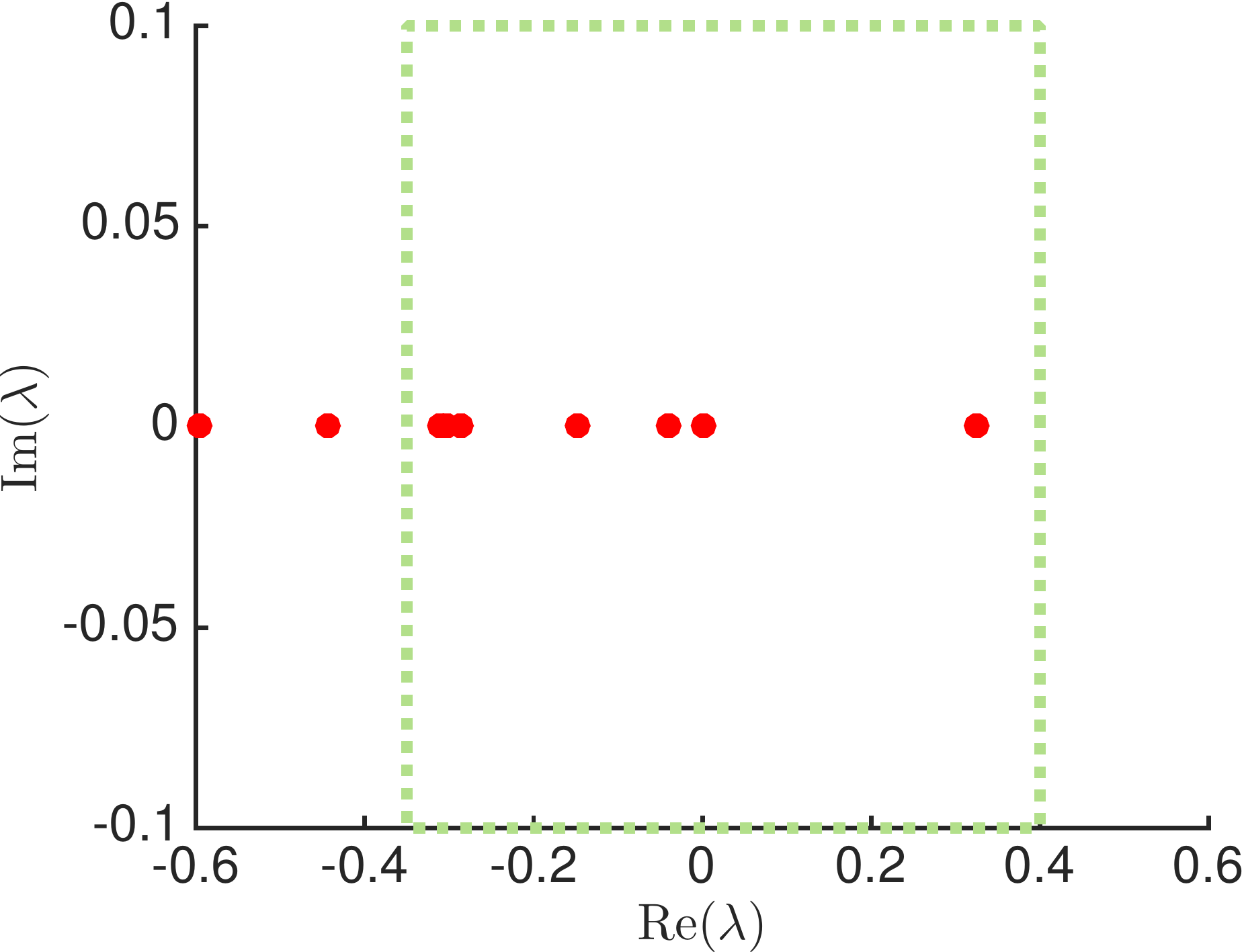} 
\includegraphics[width=0.19\textwidth]{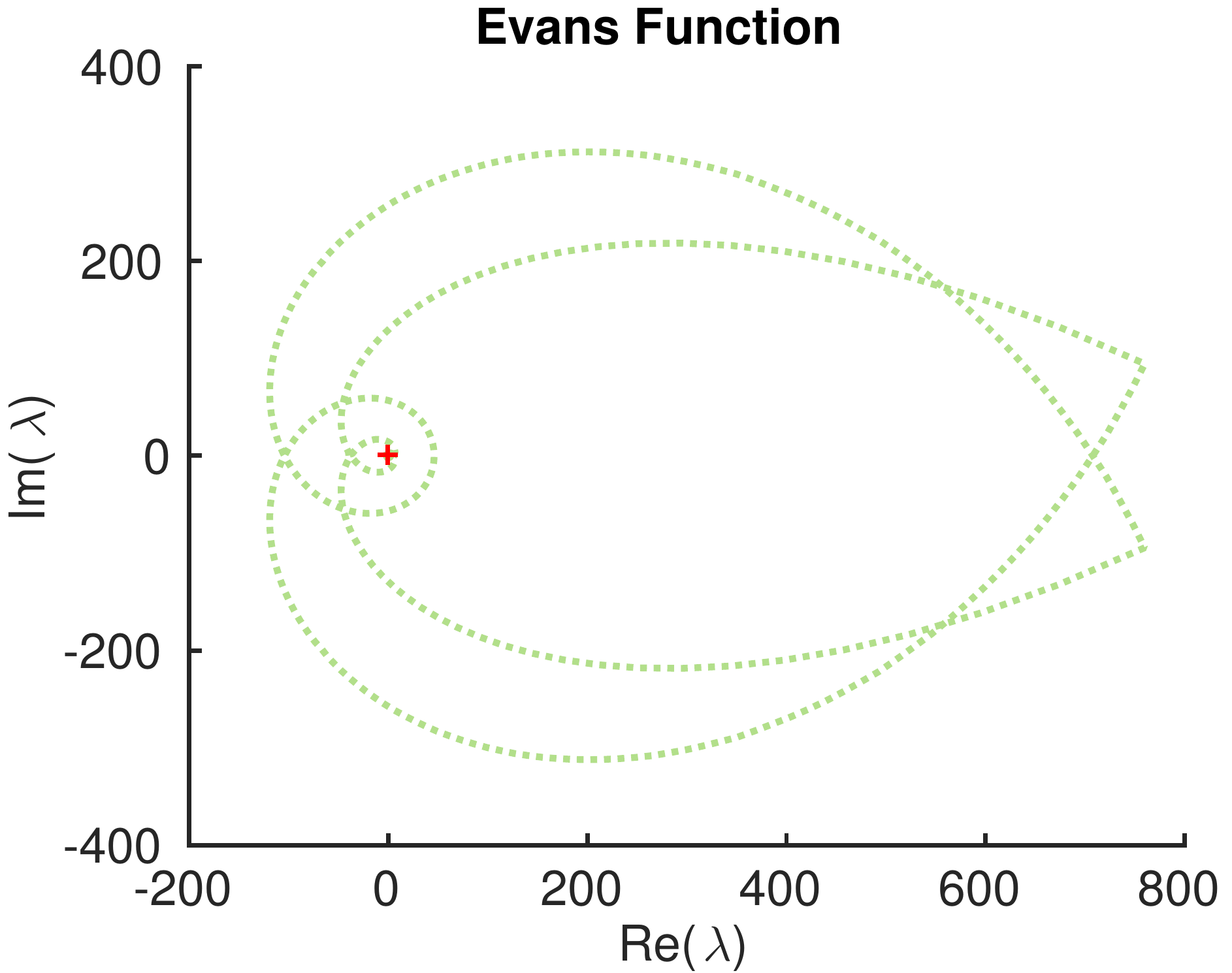} 
\includegraphics[width=0.19\textwidth]{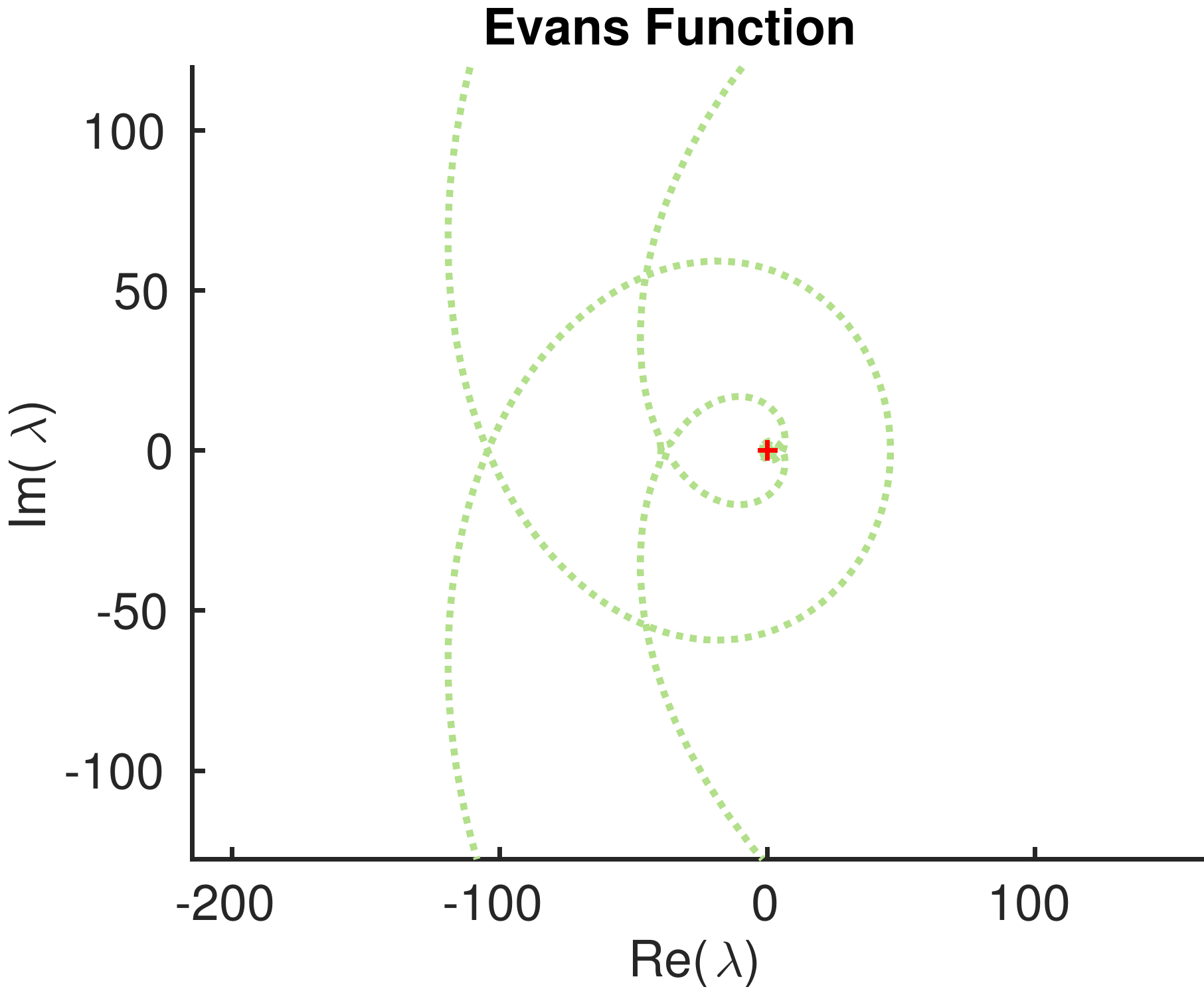}
\includegraphics[width=0.19\textwidth]{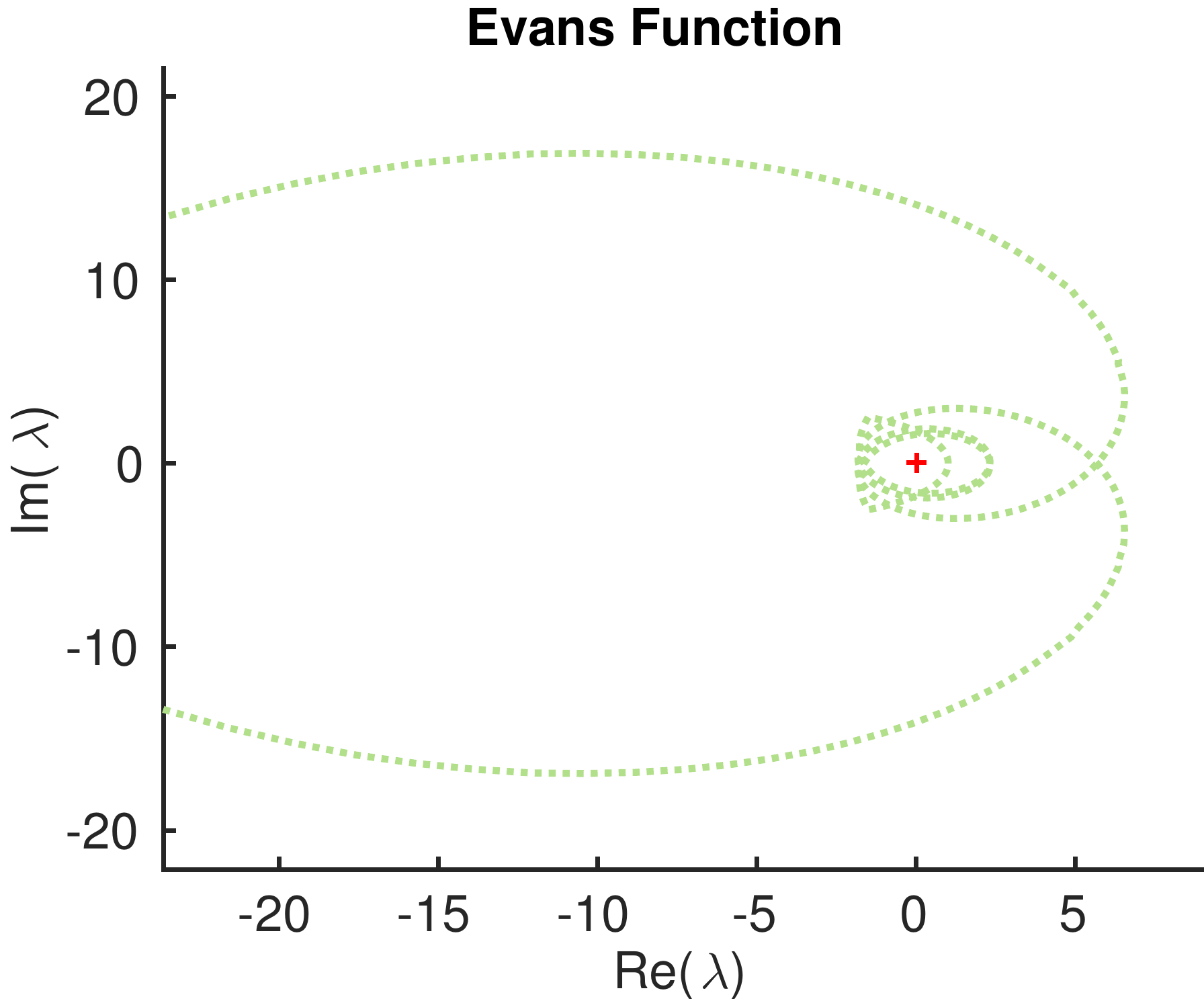}
\includegraphics[width=0.19\textwidth]{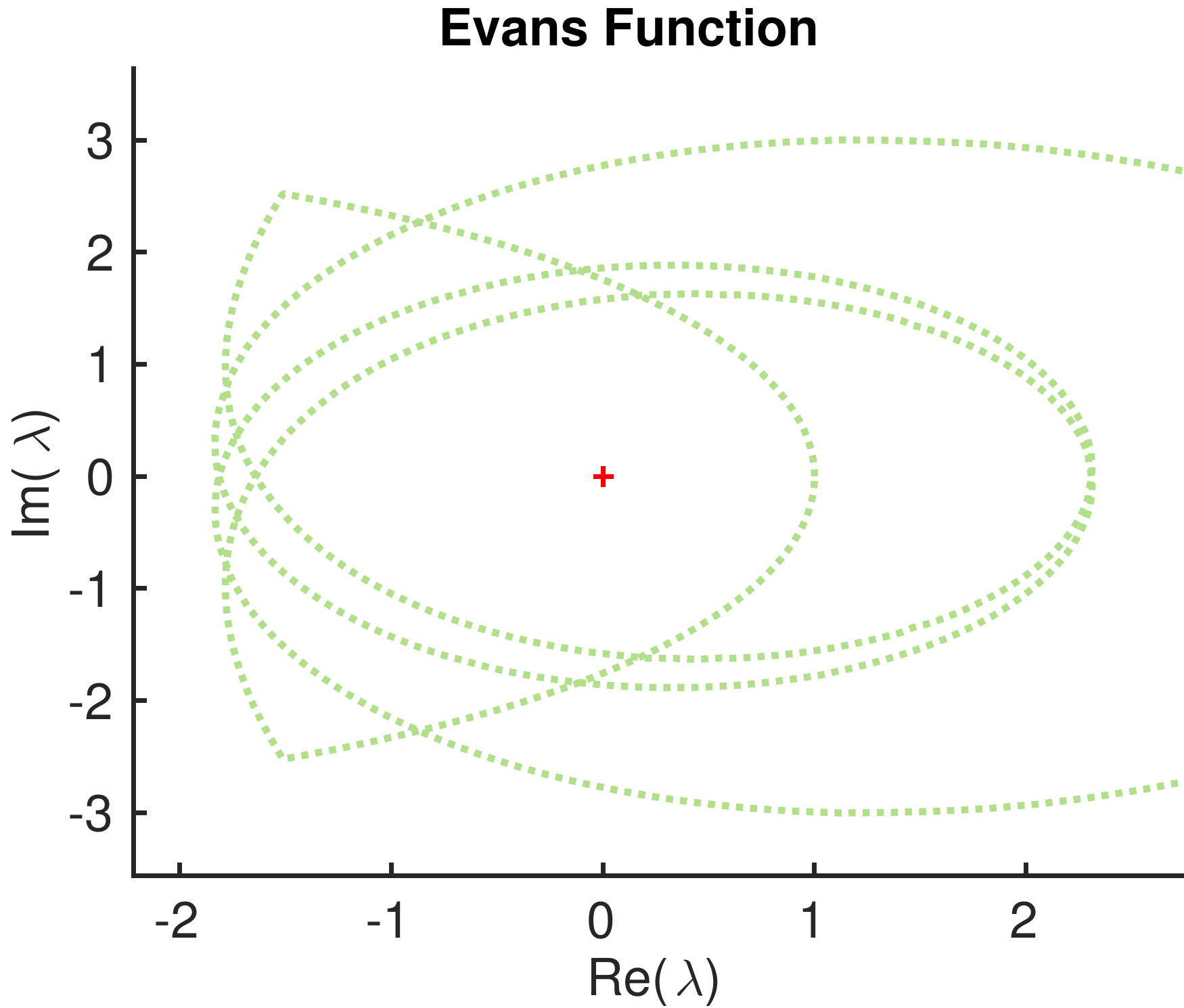} 
\caption{Shown is the Evans function of the Swift--Hohenberg equation (\ref{e:shode}) for $\mu=0.675$ and $\nu=2$. In the left panel, the computational contour is shown: it encloses 8 simple eigenvalues (indicated as red squares). The four panels to the right show all or parts of the resulting image of the contour under the Evans function with the origin indicated as a red plus sign: the resulting winding number is 8.}
\label{f:sh3}
\end{figure}

We use Fourier differentiation matrices with 8 modes to resolve the periodic $y$-direction. Linearizing the Swift--Hohenberg equation about the localized solution $u_*(x,y)$, we arrive at the eigenvalue problem
\begin{equation}\label{e:shode}
U_x = \begin{pmatrix} 0 & \id & 0 & 0 \\ 0 & 0 & \id & 0 \\
0 & 0 & 0 & \id \\ L_1 & 0 & L_2 & 0 \end{pmatrix} U
=: A(x,\lambda) U, \qquad U\in\C^{32},
\end{equation}
formulated as a first-order system in the $x$-direction, where each entry is an $8\times8$ block matrix with
\[
L_1 = -D_1 - \lambda - \mu + 3\nu u_*(x,y)^2 - 5u_*(x,y)^4, \qquad
L_2 = -2D_2,
\]
where $D_1$ and $D_2$ are the Fourier differentiation matrices approximating $(1+\Delta_y)^2$ and $1+\Delta_y$, respectively. Using an eigenvalue problem solver, we find that the Swift--Hohenberg equation with $\mu=0.675$ and $\nu=2$ linearized about the solution $u_*(x,y)$ has a double eigenvalue at $\lambda\approx0.3245$ and simple roots at
\[
\lambda\approx0,-0.0412,-0.1472,-0.2873,-0.3029,-0.3111,-0.4432,-0.5947;
\]
see \citep{Avitabile2010}. We plot the Evans function as computed with our algorithm in Figures~\ref{f:sh2} and~\ref{f:sh3}.

% --------------------------------------------------------------------------------------

\subsection{Comparison of computational performance}\label{s3.5}

It is difficult to provide a theoretical comparison of the computational costs associated with continuous orthogonalization (abbreviated by CO) and the algorithm proposed in \S\ref{s2.4} (abbreviated by BVP) as these depend strongly on a number of noncomparable parameters. Instead, we focus on comparing computational performance by using the three benchmark problems introduced in the preceding sections. We used the following procedure to compare computational times for the four implementations mentioned above separately for each of these three systems. First, we calculate the Evans function evaluated around a fixed contour to very high accuracy using continuous orthogonalization. Afterwards, we compute the Evans function around the same contour to a fixed lower accuracy using the four scripts given above and then compare the computation times for each sample system. The time comparison computations for each  system were carried out on the same computer (a MacBook Pro laptop for the coupled Nagumo system and the Korteweg--de Vries equation, and a System 76 desktop with an i7-6950x processor for the Swift--Hohenberg equation). Other computations for the Swift--Hohenberg equation were carried out with the computing cluster \textsc{oscar} at Brown University's Center for Computation and Visualization.

\begin{table}
\begin{tabular}{|c|c|c|c|c|c|c|c|}\hline
System & Parameters & L & Contour & CO & BVP (\texttt{bvp5c}) & BVP (\texttt{bvp6c}) & BVP (\texttt{bvpcheb}) \\ \hline \hline
Nagumo & $a=0.1$, $b=-1$ & 10 & $\lambda(\theta)=3+\rme^{\rmi\theta}$ & 7.94s & 40.1s & 12.4s & 3s \\ \hline
KdV & $p=4.1$, $c=5$ & 30 & $\lambda(\theta)=-0.85+\rme^{\rmi\theta}$ & 116s & 196s & 68.6s & 9.48s \\ \hline
SH & $\mu=0.675$, $\nu=2$ & 50 & $\lambda(\theta)=3+\rme^{i\theta}$ & 90.2 min & 19.0 min& 11.6 min& DNF \\
\hline
\end{tabular}
\caption{The last four columns list the times it took to compute the Evans function for the model in the first column with parameters given in the second and third columns ($L$ is the truncation length) for the contour specified in the fourth column. We did not use the solver \texttt{bvp6c} for the Swift--Hohenberg equation and also note that our solver \texttt{bvpcheb} did not work for this problem.}
\label{tb_error}
\end{table}

The comparisons are summarized in Table~\ref{tb_error}. For the coupled Nagumo system, we first computed the Evans function along the contour using the method of continuous orthogonalization with the absolute and relative error tolerances set to 1e-12. The method of continuous orthogonalization took 7.94 seconds and required an error tolerance of 1e-7 to result in a maximal relative error of 3.22e-6, while using our algorithm implemented in \texttt{bvp6c} took 12.4 seconds with an error tolerance of 1e-4 to achieve a relative error of 4.39e-6. Using \texttt{bvpcheb} required 3 seconds with degree set to 30 to obtain a relative error of 6.40e-6.

For the Korteweg--de Vries equation, we computed the Evans function on a contour of radius 1 centered at $\lambda=-0.85$. The contour had 308 points, and the relative error between successive points on the image contour was no greater than 0.1. We began by computing the Evans function using the method of continuous orthogonalization with the error tolerance in \texttt{ode15s} set to 1e-12, and took the resulting solution to represent the actual Evans function. We then computed the Evans function using continuous orthogonalization with the error tolerance set to 1e-11, our algorithm using \texttt{bvp6c} with the error tolerance set to 1e-7, and our linear BVP solver \texttt{bvpcheb} with the polynomial degree set to 60. The method of continuous orthogonalization took 116 seconds and had relative error of 8.57e-3, our algorithm using \texttt{bvp6c} took 68.6 seconds and had relative error of 7.38e-3, and our algorithm using \texttt{bvpcheb} took 9.48 seconds and had relative error of 1.55e-3.

For the Swift--Hohenberg equation on a two-dimensional domain, we computed the Evans function with an adaptive contour solver on a contour of radius 1 centered at 3. The computations required 27 points in order to achieve the relative error tolerance of 0.1 between consecutive points on the image of the contour. We found that continuous orthogonalization was fastest when the ODE tolerance was set to 1e-10. With this optimal tolerance set, it took 90.2 minutes to compute the Evans function on the contour described above. Setting the error tolerance in \texttt{bvp5c} to 1e-3, it took 19.0 minutes to compute the Evans function on the same contour. The maximum relative error between any corresponding points on the image of the Evans function for the two methods thus computed was 1.1e-2, and the average relative error was 5.6e-3. Once again, setting the error tolerance in \texttt{bvp6c} to 1e-3, it took 11.6 minutes to compute the Evans function on the same contour. The maximum relative error between any corresponding points on the image of the Evans function using bvp6c and continuous orthogonalization was 6.2e-3, and the average relative error was 2.7e-3.  While this comparison does not reflect computation times for equivalent error thresholds of the Evans function, it does indicate that our BVP algorithm may perform better (approximately 8 times faster in this example) in larger systems when achieving an acceptable error tolerance is the criteria. We also tried the linear solver \texttt{bvpcheb}, but the polynomial degree had to be large and the resulting computation time was excessive: a more sophisticated linear BVP solver might further improve computation time.

In summary, our algorithm performs well against the benchmark method of continuous orthogonalization. In particular, using a linear BVP solver has great potential with our algorithm to achieve high accuracy for low computational cost as demonstrated by the simple Chebyshev solver we built for testing purposes.

% --------------------------------------------------------------------------------------

\section{Discussion}\label{s4}

In this paper, we presented an algorithm for the computation of Evans functions that is based on solving systems of linear boundary-value problems. Our benchmark computations indicate that the algorithm performs at least as well as previous algorithms that are based on continuous orthogonalization. One of the benefits of our algorithm is that it requires only a linear boundary-value solver, which makes it easy to incorporate it into a variety of computational environments where specialized adaptive solvers may not be available.

We believe that the algorithm presented here can be used to calculate roots of Evans functions that are embedded into the essential spectrum: this is an area where eigenvalue problem solvers can generally not be used as they cannot distinguish such roots from other eigenvalues. We remark that our KdV computations show that our approach correctly captures roots of Evans functions that are extended analytically across the essential spectrum.

Another application for which we expect our algorithm to prove useful is the rigorous verification of stability of travelling waves; see \citep{Barker,SZ} for examples of recent work in this direction. For rigorous verification, our algorithm has the potentially advantageous properties that, in contrast to the exterior-product method, the dimension of the ODE system is the same as for existence, and that, in contrast to continuous orthogonalization, the ODE we need to solve is linear rather than nonlinear. Determining which of these methods works best for rigorous verification of stability is an interesting future direction.

\begin{Acknowledgment}
The computations for the Swift--Hohenberg equation were carried out using computational resources and services provided by the Center for Computation and Visualization (CCV) at Brown University: we gratefully acknowledge their support. Barker was partially supported by the National Science Foundation under grant DMS-1400872. Nguyen, Ventura, and Wahl were supported by the National Science Foundation under grant DMS-1148284. Sandstede was partially supported by the National Science Foundation under grant DMS-1408742.
\end{Acknowledgment}

% --------------------------------------------------------------------------------------	

\bibliography{references}
\bibliographystyle{sandstede}

\end{document}